\documentclass[10.5pt,a4paper]{article}

\usepackage[leqno]{amsmath}
\usepackage{amsfonts}
\usepackage{graphicx}
\usepackage{amsmath}
\usepackage{amssymb}
\usepackage{latexsym}
\usepackage{amsmath, amsfonts,amssymb, amsthm, euscript,makeidx,color,mathrsfs}
\usepackage{enumerate}

\usepackage[colorlinks,linkcolor=blue,anchorcolor=green,citecolor=red]{hyperref}

\def\sqr#1#2{{\vcenter{\vbox{\hrule height.#2pt
              \hbox{\vrule width.#2pt height#1pt \kern#1pt \vrule width.#2pt}
              \hrule height.#2pt}}}}
\def\signed #1{{\unskip\nobreak\hfil\penalty50
              \hskip2em\hbox{}\nobreak\hfil#1
              \parfillskip=0pt \finalhyphendemerits=0 \par}}
\def\endpf{\signed {$\sqr69$}}

\def\4n{\negthinspace \negthinspace \negthinspace \negthinspace }
\def\3n{\negthinspace \negthinspace \negthinspace }
\def\2n{\negthinspace \negthinspace }
\def\1n{\negthinspace }

\def\dbE{\mathbb{E}}
\def\dbF{\mathbb{F}}

\def\dbH{\mathbb{H}}

\def\dbM{\mathbb{M}}

\def\dbP{\mathbb{P}}

\def\dbR{\mathbb{R}}
\def\dbS{\mathbb{S}}

\def\dbX{\mathbb{X}}
\def\dbY{\mathbb{Y}}
\def\dbZ{\mathbb{Z}}

\def\sX{\mathscr{X}}
\def\sY{\mathscr{Y}}

\def\={\buildrel \triangle \over =}

\def\ds{\displaystyle}

\def\ns{\noalign{\ss}}
\def\a{\alpha}
\def\b{\beta}
\def\g{\gamma}
\def\d{\delta}
\def\e{\varepsilon}
\def\z{\zeta}

\def\l{\lambda}

\def\si{\sigma}

\def\f{\varphi}
\def\th{\theta}

\def\i{\infty}
\def\G{\Gamma}

\def\Th{\Theta}

\def\Si{\Sigma}
\def\F{\Phi}
\def\Om{\Omega}
\def\cA{{\cal A}}
\def\cB{{\cal B}}

\def\cD{{\cal D}}

\def\cF{{\cal F}}

\def\cH{{\cal H}}

\def\cL{{\cal L}}
\def\cM{{\cal M}}

\def\cR{{\cal R}}

\def\cU{{\cal U}}

\def\no{\noindent}

\def\ss{\smallskip}
\def\ms{\medskip}

\def\q{\quad}
\def\qq{\qquad}
\def\hb{\hbox}

\def\liminf{\mathop{\underline{\rm lim}}}

\def\Ra{\mathop{\Rightarrow}}

\def\lt{\left}
\def\rt{\right}
\def\lan{\langle}
\def\ran{\rangle}
\def\llan{\left\langle}
\def\rran{\right\rangle}
\def\blan{\big\langle}
\def\bran{\big\rangle}

\def\rf{\eqref}

\def\h{\widehat}
\def\wt{\widetilde}

\def\cd{\cdot}

\def\ae{\hbox{\rm a.e.{ }}}
\def\as{\hbox{\rm a.s.{ }}}

\def\tr{\hbox{\rm tr$\,$}}
\def\var{\hbox{\rm var$\,$}}

\def\deq{\triangleq}
\def\les{\leqslant}
\def\ges{\geqslant}

\def\({\Big (}
\def\){\Big )}
\def\[{\Big[}
\def\]{\Big]}

\def\bde{\begin{definition}\label}
\def\ede{\end{definition}}
\def\be{\begin{equation}}
\def\bel{\begin{equation}\label}
\def\ee{\end{equation}}
\def\bt{\begin{theorem}\label}
\def\et{\end{theorem}}
\def\bc{\begin{corollary}}
\def\ec{\end{corollary}}
\def\bl{\begin{lemma}\label}
\def\el{\end{lemma}}
\def\bp{\begin{proposition}\label}
\def\ep{\end{proposition}}
\def\bas{\begin{assumption}}
\def\eas{\end{assumption}}
\def\br{\begin{remark}}
\def\er{\end{remark}}
\def\bex{\begin{example}\label}
\def\ex{\end{example}}
\def\ba{\begin{array}}
\def\ea{\end{array}}
\def\ed{\end{document}}

\def\square#1{\vbox{\hrule\hbox{\vrule height#1%
     \kern#1\vrule}\hrule}}
\def\rectangle#1#2{\vbox{\hrule\hbox{\vrule height#1%
     \kern#2\vrule}\hrule}}

\font\tenbb=msbm10 \font\sevenbb=msbm7 \font\fivebb=msbm5

\newfam\bbfam
\scriptscriptfont\bbfam=\fivebb \textfont\bbfam=\tenbb
\scriptfont\bbfam=\sevenbb

\newtheorem{theorem}{\hskip 1.3em Theorem}[section]
\newtheorem{definition}[theorem]{\hskip 1.3em Definition}
\newtheorem{proposition}[theorem]{\hskip 1.3em Proposition}
\newtheorem{corollary}[theorem]{\hskip 1.3em Corollary}
\newtheorem{lemma}[theorem]{\hskip 1.3em Lemma}
\newtheorem{remark}[theorem]{\hskip 1.3em Remark}
\newtheorem{example}[theorem]{\hskip 1.3em Example}

\newtheorem{assumption}[theorem]{\hskip 1.3em Assumption}

\makeatletter
   
   \@addtoreset{equation}{section}
\makeatother

\begin{document}

\title{\bf Mean-Field Stochastic Linear Quadratic Optimal Control Problems: Open-Loop Solvabilities}
\author{Jingrui Sun\thanks{Department of Applied Mathematics, The Hong Kong Polytechnic University, Hong Kong, China
(sjr@mail.ustc.edu.cn).}}
\maketitle

\noindent {\bf Abstract:} This paper is concerned with a mean-field linear quadratic (LQ, for short)
optimal control problem with deterministic coefficients. It is shown that convexity of the cost functional
is necessary for the finiteness of the mean-field LQ problem, whereas uniform convexity of the cost functional
is sufficient for the open-loop solvability of the problem. By considering a family of uniformly convex cost
functionals, a characterization of the finiteness of the problem is derived and a minimizing sequence, whose
convergence is equivalent to the open-loop solvability of the problem, is constructed. Then, it is proved that
the uniform convexity of the cost functional is equivalent to the solvability of two coupled differential Riccati
equations and the unique open-loop optimal control admits a state feedback representation in the case that the
cost functional is uniformly convex. Finally, some examples are presented to illustrate the theory developed.


\ms


\ms

\no\bf AMS subject classifications. \rm 49N10, 49N35, 93E20

\section{Introduction}

Let $(\Om,\cF,\dbF,\dbP)$ be a complete filtered probability space on which a standard one-dimensional
Brownian motion $W=\{W(t); 0\les t<\i\}$ is defined, where $\dbF=\{\cF_t\}_{t\ges0}$ is the natural
filtration of $W$ augmented by all the $\dbP$-null sets in $\cF$.
Consider the following controlled linear stochastic differential equation (SDE, for short) on a finite horizon $[t,T]$:
\bel{state}\left\{\2n\ba{ll}
\ds dX(s)=\Big\{A(s)X(s)+\bar A(s)\dbE[X(s)]+B(s)u(s)+\bar B(s)\dbE[u(s)]+b(s)\Big\}ds\\
\ns\ds\qq\qq~+\Big\{C(s)X(s)+\bar C(s)\dbE[X(s)]+D(s)u(s)+\bar D(s)\dbE[u(s)]+\si(s)\Big\}dW(s),\qq s\in[t,T], \\
\ns\ds X(t)= \xi,\ea\right.\ee
where $A(\cd)$, $\bar A(\cd)$, $B(\cd)$, $\bar B(\cd)$, $C(\cd)$, $\bar C(\cd)$, $D(\cd)$, $\bar D(\cd)$ are given
deterministic matrix-valued functions; $b(\cd)$, $\si(\cd)$ are vector-valued $\dbF$-progressively measurable processes
and $\xi$ is an $\cF_t$-measurable random vector. In the above, $u(\cd)$ is the {\it control process} and $X(\cd)$ is the
corresponding {\it state process} with {\it initial pair} $(t,\xi)$. For any $t\in[0,T)$, we define
$$\cU[t,T]=\lt\{u:[t,T]\times\Om\to\dbR^m\bigm|u(\cd)\hb{ is $\dbF$-progressively measurable,}~\dbE\int_t^T|u(s)|^2ds<\i\rt\}.$$
Any $u(\cd)\in\cU[t,T]$ is called an {\it admissible control} (on $[t,T]$). Under some mild conditions,
for any initial pair $(t,\xi)$ with $\xi$ being square-integrable and any admissible control $u(\cd)\in\cU[t,T]$,
(\ref{state}) admits a unique square-integrable solution $X(\cd)\equiv X(\cd\,;t,\xi,u(\cd))$.
Now we introduce the following cost functional:
\bel{cost}\ba{ll}
\ds J(t,\xi;u(\cd))\deq \dbE\Bigg\{\lan GX(T),X(T)\ran+2\lan g,X(T)\ran+\lan\bar G\dbE[X(T)],\dbE[X(T)]\ran+2\lan\bar g,\dbE[X(T)]\ran\\
\ns\ds\qq\qq\,+\int_t^T\left[\llan\begin{pmatrix}Q(s)&S(s)^\top\\S(s)&R(s)\end{pmatrix}
                                \begin{pmatrix}X(s)\\ u(s)\end{pmatrix},
                                \begin{pmatrix}X(s)\\u(s)\end{pmatrix}\rran
+2\llan\begin{pmatrix}q(s)\\ \rho(s)\end{pmatrix},\begin{pmatrix}X(s)\\ u(s)\end{pmatrix}\rran\right]ds\\
\ns\ds\qq\qq\,+\int_t^T\left[\llan\begin{pmatrix}\bar Q(s)&\bar S(s)^\top\\\bar S(s)&\bar R(s)\end{pmatrix}
                                \begin{pmatrix}\dbE[X(s)]\\ \dbE[u(s)]\end{pmatrix},
                                \begin{pmatrix}\dbE[X(s)]\\\dbE[u(s)]\end{pmatrix}\rran
+2\llan\begin{pmatrix}\bar q(s)\\ \bar\rho(s)\end{pmatrix},\begin{pmatrix}\dbE[X(s)]\\\dbE[u(s)]\end{pmatrix}\rran\right] ds\Bigg\},
\ea\ee
where $G$, $\bar G$ are symmetric matrices and $Q(\cd)$, $\bar Q(\cd)$, $S(\cd)$, $\bar S(\cd)$, $R(\cd)$, $\bar R(\cd)$ are
deterministic matrix-valued functions with $Q(\cd)^\top=Q(\cd)$, $\bar Q(\cd)^\top=\bar Q(\cd)$,
$R(\cd)^\top=R(\cd)$, $\bar R(\cd)^\top=\bar R(\cd)$; $g$ is an $\cF_T$-measurable random vector and $\bar g$ is a (deterministic) vector;
$q(\cd)$, $\rho(\cd)$ are vector-valued $\dbF$-progressively measurable processes and $\bar q(\cd)$, $\bar \rho(\cd)$
are deterministic vector-valued functions. Our mean-field stochastic LQ optimal control problem can be stated as follows:

\ms

\bf Problem (MF-LQ). \rm For any given initial pair $(t,\xi)\in[0,T)\times L^2_{\cF_t}(\Om;\dbR^n)$, find a $u^*(\cd)\in\cU[t,T]$ such that
\bel{optim}J(t,\xi;u^*(\cd))=\inf_{u(\cd)\in\cU[t,T]}J(t,\xi;u(\cd))\deq V(t,\xi).\ee

In the above, $L^2_{\cF_t}(\Om;\dbR^n)$ is the space of all $\cF_t$-measurable, $\dbR^n$-valued random vectors $\xi$
with $\dbE|\xi|^2<\i$. Any $u^*(\cd)\in\cU[t,T]$ satisfying (\ref{optim}) is called an ({\it open-loop}) {\it optimal control} of
Problem (MF-LQ) for the initial pair $(t,\xi)$, and the corresponding $X^*(\cd)\equiv X(\cd\,; t,\xi,u^*(\cd))$
is called an {\it optimal state process}. The function $V(\cd\,,\cd)$ is called the {\it value function} of Problem (MF-LQ).
In the special case of $b(\cd)$, $\si(\cd)$, $g(\cd)$, $\bar g(\cd)$, $q(\cd)$, $\bar q(\cd)$, $\rho(\cd)$, $\bar \rho(\cd)=0$,
we denote by $J^0(t,\xi; u(\cd))$, $V^0(t,\xi)$ and Problem (MF-LQ)$^0$ the corresponding cost functional, value function and Problem (MF-LQ),
respectively.

\ms

Comparing with the classical stochastic LQ optimal control problem, a new feature of Problem (MF-LQ) is that both the state equation
and the cost functional involve the states and the controls as well as their expectations. In this case,
we call \rf{state} a controlled mean-field (forward) SDE (MF-SDE, for short).
The history of MF-SDEs can be traced back to the work of Kac \cite {Kac 1956} in 1956 and McKean \cite{McKean 1966} in 1966.
Since then, many researchers have made contributions to such kind of equations and applications; see, for example, Dawson \cite{Dawson 1983},
Dawson--G$\ddot{{\rm a}}$rtner \cite{Dawson-Gartner 1987}, Scheutzow \cite{Scheutzow 1987}, G$\ddot{{\rm a}}$rtner \cite{Gartner 1988},
Graham \cite{Graham 1992}, Chan \cite{Chan 1994}, Chiang \cite{Chiang 1994} and Ahmed--Ding \cite{Ahmed-Ding 1995}.
For recent development of MF-SDEs, readers may refer to Huang--Malham\'{e}--Caines \cite{Huang-Malhame-Caines 2006},
Veretennikov \cite{Veretennikov 2006}, Mahmudov--McKibben \cite{Mahmudov-McKibben 2007},
Buckdahn--Djehiche--Li--Peng \cite{Buckdahn-Djehiche-Li-Peng 2009}, Buckdahn--Li--Peng \cite{Buckdahn-Li-Peng 2009},
Borkar--Kumar \cite{Borkar-Kumar 2010}, Crisan--Xiong \cite{Crisan-Xiong 2010}, Kotelenez--Kurtz \cite{Kotelenez-Kurtz 2010}
and the references cited therein.
Control problems of MF-SDEs were studied by Ahmed--Ding \cite{Ahmed-Ding 2001}, Ahmed \cite{Ahmed 2007},
Park--Balasubramaniam--Kang \cite{Park-Balasubramaniam-Kang 2008}, Buckdahn--Djehiche--Li \cite{Buckdahn-Djehiche-Li 2011},
Andersson--Djehiche \cite{Andersson-Djehiche 2011}, Meyer-Brandis--${\O}$ksendal--Zhou \cite{Meyer-Oksendal-Zhou 2012}, and so on.
More recently, Yong \cite{Yong 2013} investigated an LQ problem for MF-SDEs in finite horizons and gave some interesting motivation
for the control problem with $\dbE[X(\cd)]$ and $\dbE[u(\cd)]$ being included in the cost functional. Later, Huang--Li--Yong
\cite{Huang-Li-Yong 2014} generalized the results in \cite{Yong 2013} to the case with an infinite time horizon.

\ms

In \cite{Yong 2013}, two coupled differential Riccati equations are derived by decoupling the optimality system.
It is shown that under certain conditions, the two Riccati equations are uniquely solvable and Problem (MF-LQ)
admits a unique optimal control which has a state feedback representation. To be precise, if
\bel{Classic}\left\{\2n\ba{ll}
\ds G,~G+\bar G\ges0,\q &\ds Q(s),~Q(s)+\bar Q(s)\ges0,\\
\ns\ds S(s)=\bar S(s)=0,\q &\ds R(s),~R(s)+\bar R(s)\ges\d I,
\ea\right.\qq\ae~s\in[0,T],\ee
for some $\d>0$, then the unique solvability of the two Riccati equations can be obtained from the classical result
\cite[Theorem 7.2]{Yong-Zhou 1999}. However, examples show that the two Riccati equations might still be solvable even
if both $R(\cd)$ and $\bar R(\cd)$ are negative semi-definite (see Example \ref{exm1}). On the other hand, it may happen
that Problem (MF-LQ) is open-loop solvable, while the optimal control cannot be obtained by solving the corresponding Riccati
equations due to the possible singularities of the terms $R+D^\top PD$ and $R+\bar R+(D+\bar D)^\top P(D+\bar D)$
(see Example \ref{ex2}). Thus, some questions arise naturally:
(a) What is the relationship between Problem (MF-LQ) and the solvability of the two Riccati equations?
(b) How can we characterize the open-loop solvability of Problem (MF-LQ)?
(c) How can we find an optimal control in general?
The purpose of this paper is to study Problem (MF-LQ) from an open-loop point of view and to address the above issues.
Closed-loop mean-field LQ problems will be investigated in a forthcoming paper.

\ms

Our main idea and results of this paper can be informally described as follows. By a representation of the cost functional,
we first show that for the open-loop solvability of Problem (MF-LQ), a necessary condition is the convexity of the cost functional
and a sufficient condition is the uniform convexity of the cost functional. Under the convexity condition, by adding
$\e\dbE\int_t^T|u(s)|^2ds$ ($\e>0$) to the original cost functional, we get a family of uniformly convex functionals. The corresponding
mean-field LQ problems admit unique optimal controls $u^*_\e(\cd),~\e>0$, which form a minimizing sequence of Problem (MF-LQ). Then the
open-loop solvability of Problem (MF-LQ) is characterized by the convergence of the sequence, whose limit is an optimal control of
Problem (MF-LQ). To construct $u^*_\e(\cd)$ explicitly, we further investigate Problem (MF-LQ) with uniformly convex cost functionals.
Since the uniform convexity condition is much weaker than \rf{Classic}, the result in \cite{Yong 2013} fails to apply to this case.
To overcome this difficulty, we reduce Problem (MF-LQ) to a classical stochastic LQ problem and a deterministic LQ problem.
By making use of a result found in \cite{Sun-Li-Yong}, we establish the equivalence between the uniform convexity of the cost functional
and the solvability of the two Riccati equations. Then by the completion of squares technique, we obtain a state feedback representation
of the optimal control via the solutions of the two Riccati equations.

\ms

The rest of the paper is organized as follows. Section 2 gives some preliminaries. In Section 3,
we study Problem (MF-LQ) from a Hilbert space viewpoint and derive necessary and sufficient conditions
for the finiteness and open-loop solvability of the problem by considering a family of uniformly convex cost functionals.
Section 4 shows that the solvability of two coupled Riccati equations is necessary for the uniform convexity of
the cost functional. In Section 5, we further prove that the solvability of the two coupled Riccati equations is
also sufficient for the uniform convexity of the cost functional. Moreover, a state feedback representation is obtained for
the optimal control. Some illustrative examples are presented in Section 6.

\section{Preliminaries}

Throughout this paper, we denote by $\dbR^{n\times m}$ the Euclidean space of all $n\times m$ real matrices,
and by $\dbS^n$ the space of all symmetric $n\times n$ real matrices. Recall that the inner product
$\lan\cd\,,\cd\ran$ on $\dbR^{n\times m}$ is given by $\lan M,N\ran\mapsto\tr(M^\top N)$, where the superscript
$\top$ denotes the transpose of vectors or matrices, and the induced norm is given by $|M|=\sqrt{\tr(M^\top M)}$.
When there is no confusion, we shall use $\lan\cd\,,\cd\ran$ for inner products in possibly different Hilbert spaces,
and denote by $|\cd|$ the norm induced by $\lan\cd\,,\cd\ran$.
For a matrix $M\in\dbR^{n\times m}$, we denote by $\cR(M)$ the range of $M$, and if $M\in\dbS^n$, we use the notation
$M>0~(\ges0)$ to indicate that $M$ is positive (semi-) definite. For a bounded linear operator $\cA$ form a Banach
$\sX$ space into another Banach space $\sY$, we denote by $\cA^*$ the adjoint operator of $\cA$.
Let $T>0$ be a fixed time horizon. For any $t\in[0,T]$
and Euclidean space $\dbH$, we let $L^p(t,T;\dbH)$ $(1\les p\les\i)$ be the space of all $\dbH$-valued functions
that are $L^p$-integrable on $[t,T]$ and $C([t,T];\dbH)$ be the space of all $\dbH$-valued continuous functions on $[t,T]$.
Next, we introduce the following spaces:
$$\ba{ll}
\ds L^2_{\cF_t}(\Om;\dbH)=\Big\{\xi:\Om\to\dbH\bigm|\xi\hb{ is $\cF_t$-measurable, }\dbE|\xi|^2<\i\Big\},\\
\ns\ds L_\dbF^2(t,T;\dbH)=\lt\{\f:[t,T]\times\Om\to\dbH\bigm|\f(\cd)\hb{ is
$\dbF$-progressively measurable, }\dbE\int^T_t|\f(s)|^2ds<\i\rt\},\\
\ns\ds L_\dbF^2(\Om;C([t,T];\dbH))=\lt\{\f:[t,T]\times\Om\to\dbH\bigm|\f(\cd)\hb{
is $\dbF$-adapted, continuous, }\dbE\lt(\sup_{s\in[t,T]}|\f(s)|^2\rt)<\i\rt\},\\
\ns\ds L^2_\dbF(\Om;L^1(t,T;\dbH))=\lt\{\f:[t,T]\times\Om\to\dbH\bigm|\f(\cd)\hb{ is $\dbF$-progressively measurable, }
\dbE\lt(\int_t^T|\f(s)|ds\rt)^2<\i\rt\}.\ea$$
Further, we introduce the following notation: For any $\dbS^n$-valued measurable function $F$ on $[t,T]$,
$$\left\{\2n\ba{ll}
\ds F\ges0\q\Longleftrightarrow\q F(s)\ges0,\qq\ae~s\in[t,T],\\
\ns\ds F>0\q\Longleftrightarrow\q F(s)>0,\qq\ae~s\in[t,T],\\
\ns\ds F\gg0\q\Longleftrightarrow\q F(s)\ges\d I,\qq\ae~s\in[t,T],\hb{ for some }\d>0.
\ea\right.$$

\ms

The following assumptions will be in force throughout this paper.

\ms

{\bf(H1)} The coefficients of the state equation satisfy the following:
$$\left\{\2n\ba{ll}
\ds A(\cd),\bar{A}(\cd)\in L^1(0,T;\dbR^{n\times n}),\q B(\cd),\bar{B}(\cd)\in L^2(0,T;\dbR^{n\times m}),
\q b(\cd)\in L^2_\dbF(\Om;L^1(0,T;\dbR^n)),\\
\ns\ds C(\cd),\bar C(\cd)\in L^2(0,T;\dbR^{n\times n}),\q D(\cd),\bar D(\cd)\in L^{\i}(0,T;\dbR^{n\times m}),
\q\si(\cd)\in L^2_\dbF(0,T;\dbR^n).
\ea\right.$$

{\bf(H2)} The weighting coefficients in the cost functional satisfy the following:
$$\left\{\2n\ba{ll}
\ds Q(\cd), \bar Q(\cd)\in L^1(0,T;\dbS^n),\q S(\cd), \bar S(\cd)\in L^2(0,T;\dbR^{m\times n}),
\q R(\cd), \bar R(\cd)\in L^\i(0,T;\dbS^m),\\
\ns\ds g\in L^2_{\cF_T}(\Om;\dbR^n),\q q(\cd)\in L^2_\dbF(\Om;L^1(0,T;\dbR^n)),\q \rho(\cd)\in L_\dbF^2(0,T;\dbR^m),\\
\ns\ds\bar g\in \dbR^n,\q \bar q(\cd)\in L^1(0,T;\dbR^n),\q \bar\rho(\cd)\in L^2(0,T;\dbR^m),\q G, \bar G\in\dbS^n.
\ea\right.$$

\ms

By a standard argument using contraction mapping theorem, one can show that under (H1),
for any $(t,\xi)\in[0,T)\times L^2_{\cF_t}(\Om;\dbR^n)$ and any $u(\cd)\in\cU[t,T]$,
(\ref{state}) admits a unique solution $X(\cd)\equiv X(\cd\,; t,\xi,u(\cd))\in L_\dbF^2(\Om;C([t,T];\dbR^n))$.
Hence, under (H1)--(H2), the cost functional (\ref{cost}) is well-defined, and Problem (MF-LQ) makes sense.
Now we introduce the following definition.

\bde{def-1} \rm (i) Problem (MF-LQ) is said to be {\it finite at initial pair $(t,\xi)\in[0,T]\times L^2_{\cF_t}(\Om;\dbR^n)$} if
\bel{def-1-1}V(t,\xi)>-\i.\ee
Problem (MF-LQ) is said to be {\it finite at $t\in[0,T]$} if \rf{def-1-1} holds for all $\xi\in L^2_{\cF_t}(\Om;\dbR^n)$,
and Problem (MF-LQ) is said to be {\it finite} if it is finite at all $t\in[0,T]$.

\ms

(ii) Problem (MF-LQ) is said to be ({\it uniquely}) {\it open-loop solvable at initial pair $(t,\xi)\in[0,T]\times L^2_{\cF_t}(\Om;\dbR^n)$}
if there exists a (unique) $u^*(\cd)\in\cU[t,T]$ satisfying \rf{optim}. Problem (MF-LQ) is said to be ({\it uniquely})
{\it open-loop solvable at $t$} if for any $\xi\in L^2_{\cF_t}(\Om;\dbR^n)$, there exists a (unique) $u^*(\cd)\in\cU[t,T]$ satisfying \rf{optim},
and Problem (MF-LQ) is said to be ({\it uniquely}) {\it open-loop solvable} ({\it on $[0,T)$}) if it is (uniquely) open-loop solvable at all $t\in[0,T)$.
\ede

Next, we introduce the following mean-field backward SDE (MF-BSDE, for short) associated with the state process
$X(\cd)\equiv X(\cd\,; t,\xi,u(\cd))$:
\bel{adj-X}\left\{\2n\ba{ll}
\ds dY(s)=-\Big\{A^\top Y+\bar A^\top\dbE[Y]+C^\top Z+\bar C^\top\dbE[Z]+QX+\bar Q\dbE[X]\\
\ns\ds\qq\qq\q~~+S^\top u+\bar S^\top\dbE[u]+q+\bar q\Big\}ds+ZdW(s), \qq s\in[t,T],\\
\ns\ds Y(T)=GX(T)+\bar G\dbE[X(T)]+g+\bar g.
\ea\right.\ee
The following result is concerned with the differentiability of the map $u(\cd)\mapsto J(t,\xi;u(\cd))$.

\bp{u+v} \sl Let {\rm(H1)--(H2)} hold and $t\in[0,T)$ be given.
For any $\xi\in L^2_{\cF_t}(\Om;\dbR^n), \l\in\dbR$ and $u(\cd), v(\cd)\in\cU[t,T]$,
the following holds:
\bel{u+lamv}\ba{ll}
\ds J(t,\xi;u(\cd)+\l v(\cd))-J(t,\xi;u(\cd))\\
\ns\ds=\l^2 J^0(t,0;v(\cd))+2\l\dbE\int_t^T\blan B^\top Y+\bar B^\top\dbE[Y]+D^\top Z+\bar D^\top\dbE[Z]\\
\ns\ds\qq\qq\qq\qq\qq\qq\q~+SX+\bar S\dbE[X]+Ru+\bar R\dbE[u]+\rho+\bar\rho,v\bran ds,
\ea\ee
where $X(\cd)=X(\cd\,;t,\xi,u(\cd))$ and $(Y(\cd),Z(\cd))$ is the adapted solution to the MF-BSDE \rf{adj-X} associated with $X(\cd)$.
Consequently, the map $u(\cd)\mapsto J(t,\xi;u(\cd))$ is Fr\'echet differentiable with the Fr\'echet derivative given by
\bel{DJ}\ba{ll}
\ds\cD J(t,\xi;u(\cd))(s)=2\Big\{B(s)^\top Y(s)+\bar B(s)^\top\dbE[Y(s)]+D(s)^\top Z(s)+\bar D(s)^\top\dbE[Z(s)]+S(s)X(s)\\
\ns\ds\qq\qq\qq\qq\q~
+\bar S(s)\dbE[X(s)]+R(s)u(s)+\bar R(s)\dbE[u(s)]+\rho(s)+\bar\rho(s)\big]\Big\},\qq s\in[t,T].\ea\ee
\ep

\it Proof. \rm Let $\h X(\cd)=X(\cd\,;t,\xi,u(\cd)+\l v(\cd))$ and $X_0(\cd)$ be the solution to the following MF-SDE:
$$\left\{\2n\ba{ll}
\ds dX_0(s)=\Big\{A(s)X_0(s)+\bar A(s)\dbE[X_0(s)]+B(s)v(s)+\bar B(s)\dbE[v(s)]\Big\}ds\\
\ns\ds\qq\qq~+\Big\{C(s)X_0(s)+\bar C(s)\dbE[X_0(s)]+D(s)v(s)+\bar D(s)\dbE[v(s)]\Big\}dW(s),\qq s\in[t,T], \\
\ns\ds X_0(t)=0,\ea\right.$$
By the linearity of the state equation, $\h X(\cd)=X(\cd)+\l X_0(\cd)$. Hence,
$$\ba{ll}
\ds J(t,\xi;u(\cd)+\l v(\cd))-J(t,\xi;u(\cd))\\
\ns\ds=\l\dbE\Bigg\{\blan G\big[2X(T)+\l X_0(T)\big],X_0(T)\bran+2\lan g,X_0(T)\ran\\
\ns\ds\qq\q~~+\int_t^T\left[\llan\begin{pmatrix}Q&S^\top\\S&R\end{pmatrix}
                              \begin{pmatrix}2X+\l X_0\\2u+\l v\end{pmatrix},
                              \begin{pmatrix}X_0\\v\end{pmatrix}\rran
+2\llan\begin{pmatrix}q\\\rho\end{pmatrix},\begin{pmatrix}X_0\\v\end{pmatrix}\rran\right]ds\Bigg\}\\
\ns\ds~~+\l\Bigg\{\llan\bar G\(2\dbE[X(T)]+\l\dbE[X_0(T)]\),\dbE[X_0(T)]\rran+2\lan\bar g,\dbE[X_0(T)]\ran\\
\ns\ds\qq\q~~+\int_t^T\left[\llan\begin{pmatrix}\bar Q&\bar S^\top\\\bar S&\bar R\end{pmatrix}
                              \begin{pmatrix}2\dbE[X]+\l \dbE[X_0]\\2\dbE[u]+\l \dbE[v]\end{pmatrix},
                              \begin{pmatrix}\dbE[X_0]\\\dbE[v]\end{pmatrix}\rran
+2\llan\begin{pmatrix}\bar q\\\bar\rho\end{pmatrix},\begin{pmatrix}\dbE[X_0]\\\dbE[v]\end{pmatrix}\rran\right]ds\Bigg\}\\
\ns\ds=2\l\dbE\Bigg\{\lan GX(T)+g,X_0(T)\ran+\int_t^T\[\lan QX+S^\top u+q,X_0\ran+\lan SX+Ru+\rho,v\ran\]ds\Bigg\}\\
\ns\ds~~+\l^2\dbE\Bigg\{\lan GX_0(T),X_0(T)\ran+\int_t^T\llan\begin{pmatrix}Q&S^\top\\S&R\end{pmatrix}
                                                        \begin{pmatrix}X_0\\ v\end{pmatrix},
                                                        \begin{pmatrix}X_0\\ v\end{pmatrix}\rran ds\Bigg\}
\ea$$
$$\ba{ll}
\ds~~+2\l\Bigg\{\lan\bar G\dbE[X(T)]\1n+\1n\bar g,\dbE[X_0(T)]\ran
\1n+\1n\int_t^T\[\lan\bar Q\dbE[X]\1n+\1n\bar S^\top\dbE[u]\1n+\1n\bar q,\dbE[X_0]\ran\1n+\1n\lan\bar S\dbE[X]\1n+\1n\bar R\dbE[u]\1n+\1n\bar\rho,\dbE[v]\ran\]ds\Bigg\}\\
\ns\ds~~+\l^2\Bigg\{\lan\bar G\dbE[X_0(T)],\dbE[X_0(T)]\ran
+\int_t^T\llan\begin{pmatrix}\bar Q&\bar S^\top\\\bar S&\bar R\end{pmatrix}
             \begin{pmatrix}\dbE[X_0]\\ \dbE[v]\end{pmatrix},
             \begin{pmatrix}\dbE[X_0]\\ \dbE[v]\end{pmatrix}\rran ds\Bigg\}\\
\ns\ds=2\l\dbE\Bigg\{\lan GX(T)+\bar G\dbE[X(T)]+g+\bar g,X_0(T)\ran\\
\ns\ds\qq\q~~+\int_t^T\[\blan QX\1n+\1n\bar Q\dbE[X]\1n+\1nS^\top u\1n+\1n\bar S^\top\dbE[u]\1n+\1nq\1n+\1n\bar q,X_0\bran
\1n+\1n\blan SX\1n+\1n\bar S\dbE[X]\1n+\1nRu\1n+\1n\bar R\dbE[u]\1n+\1n\rho\1n+\1n\bar\rho,v\bran\]ds\Bigg\}\\
\ns\ds~~+\l^2 J^0(t,0;v(\cd)).
\ea$$
Now applying It\^o's formula to $s\mapsto\lan Y(s),X_0(s)\ran$, we have
$$\ba{ll}
\ds\dbE\lan GX(T)+\bar G\dbE[X(T)]+g+\bar g,X_0(T)\ran\\
\ns\ds=\dbE\int_t^T\Big\{-\lan A^\top Y+\bar A^\top\dbE[Y]+C^\top Z+\bar C^\top\dbE[Z]+QX+\bar Q\dbE[X]
+S^\top u+\bar S^\top\dbE[u]+q+\bar q,X_0\ran\\
\ns\ds\qq\qq~+\lan AX_0+\bar A\dbE[X_0]+Bv+\bar B\dbE[v],Y\ran+\lan CX_0+\bar C\dbE[X_0]+Dv+\bar D\dbE[v],Z\ran\Big\}ds\\
\ns\ds=\dbE\int_t^T\Big\{\lan B^\top Y+\bar B^\top\dbE[Y]+D^\top Z+\bar D^\top\dbE[Z],v\ran
-\lan QX+\bar Q\dbE[X]+S^\top u+\bar S^\top\dbE[u]+q+\bar q,X_0\ran\Big\}ds.
\ea$$
Combining the above equalities, we obtain \rf{u+lamv}.
\endpf

\ms

From the above, we have the following result, which gives a characterization for the optimal controls of Problem (MF-LQ).

\bt{}\sl Let {\rm(H1)--(H2)} hold and $(t,\xi)\in[0,T)\times L^2_{\cF_t}(\Om;\dbR^n)$ be given. Let
$u^*(\cd)\in\cU[t,T]$ and $(X^*(\cd),Y^*(\cd),Z^*(\cd))$ be the adapted solution to the following
(decoupled) mean-field forward-backward stochastic differential equation (MF-FBSDE, for short):
\bel{9-14-MFBSDE}\left\{\2n\ba{ll}
\ds dX^*(s)=\Big\{AX^*+\bar A\dbE[X^*]+Bu^*+\bar B\dbE[u^*]+b\Big\}ds\\
\ns\ds\qq\qq~+\Big\{CX^*+\bar C\dbE[X^*]+Du^*+\bar D\dbE[u^*]+\si\Big\}dW(s),\qq s\in[t,T], \\
\ns\ds dY^*(s)=-\Big\{A^\top Y^*+\bar A^\top\dbE[Y^*]+C^\top Z^*+\bar C^\top\dbE[Z^*]+QX^*+\bar Q\dbE[X^*]\\
\ns\ds\qq\qq\q~~+S^\top u^*+\bar S^\top\dbE[u^*]+q+\bar q\Big\}ds+Z^*dW(s), \qq s\in[t,T],\\
\ns\ds X^*(t)=\xi,\qq Y^*(T)=GX^*(T)+\bar G\dbE[X^*(T)]+g+\bar g.
\ea\right.\ee
Then $u^*(\cd)$ is an optimal control of Problem {\rm(MF-LQ)} for the initial pair $(t,\xi)$ if and only if
\bel{9-14-tu}J^0(t,0;u(\cd))\ges0,\qq\forall u(\cd)\in\cU[t,T],\ee
and the following stationarity condition holds:
\bel{9-14-pingwen}\ba{ll}
\ds \cD J(t,\xi;u^*(\cd))=2\Big\{B^\top Y^*+D^\top Z^*+SX^*+ Ru^*+\rho\\
\ns\ds\qq\qq\qq\qq~+\bar B^\top\dbE[Y^*]+\bar D^\top\dbE[Z^*]+\bar S\dbE[X^*]+\bar R\dbE[u^*]+\bar\rho\Big\}=0,\q\ae~\as
\ea\ee
\et

\it Proof. \rm By \rf{u+lamv}, we see that $u^*(\cd)$ is an optimal control of Problem {\rm(MF-LQ)}
for the initial pair $(t,\xi)$ if and only if
$$\ba{ll}
\ds\l^2 J^0(t,0;u(\cd))+\l\dbE\int_t^T\lan \cD J(t,\xi;u^*(\cd))(s),u(s)\ran ds\\
\ns\ds=J(t,\xi;u^*(\cd)+\l u(\cd))-J(t,\xi;u^*(\cd))\ges0,\qq\forall \l\in\dbR,\q\forall u(\cd)\in\cU[t,T],
\ea$$
which is equivalent to \rf{9-14-tu} and the following:
$$\dbE\int_t^T\lan \cD J(t,\xi;u^*(\cd))(s),u(s)\ran ds\les0,\qq \forall u(\cd)\in\cU[t,T].$$
Note that the above inequality holds for all $u(\cd)\in\cU[t,T]$ if and only if $\cD J(t,\xi;u^*(\cd))(\cd)=0$.
The result therefore follows.
\endpf

\section{Finiteness and Open-Loop Solvability of Problem (MF-LQ)}

We begin with a representation of the cost functional. For any $u(\cd)\in\cU[t,T]$,
let $X_0^u(\cd)$ be the solution of
\bel{X0u}\left\{\2n\ba{ll}
\ds dX_0^u(s)=\Big\{A(s)X_0^u(s)+\bar A(s)\dbE[X_0^u(s)]+B(s)u(s)+\bar B(s)\dbE[u(s)]\Big\}ds\\
\ns\ds\qq\qq~~+\Big\{C(s)X_0^u(s)+\bar C(s)\dbE[X_0^u(s)]+D(s)u(s)+\bar D(s)\dbE[u(s)]\Big\}dW(s),\qq s\in[t,T], \\
\ns\ds X_0^u(t)=0.\ea\right.\ee
By the linearity of \rf{X0u}, we can define bounded linear operators
$\cL_t:\cU[t,T]\to L_\dbF^2(t,T;\dbR^n)$ and $\h\cL_t:\cU[t,T]\to L^2_{\cF_T}(\Om;\dbR^n)$
by $u(\cd)\mapsto X_0^u(\cd)$ and $u(\cd)\mapsto X_0^u(T)$, respectively, via the MF-SDE \rf{X0u}.
Then
$$\ba{ll}
\ds J^0(t,0;u(\cd))=\dbE\Bigg\{\lan GX_0^u(T),X_0^u(T)\ran+\lan\bar G\dbE[X_0^u(T)],\dbE[X_0^u(T)]\ran\\
\ns\ds\qq\qq\qq\qq\,+\int_t^T\llan\begin{pmatrix}Q(s)&S(s)^\top\\S(s)&R(s)\end{pmatrix}
                                \begin{pmatrix}X_0^u(s)\\ u(s)\end{pmatrix},
                                \begin{pmatrix}X_0^u(s)\\u(s)\end{pmatrix}\rran ds\\
\ns\ds\qq\qq\qq\qq\,+\int_t^T\llan\begin{pmatrix}\bar Q(s)&\bar S(s)^\top\\\bar S(s)&\bar R(s)\end{pmatrix}
                                \begin{pmatrix}\dbE[X_0^u(s)]\\ \dbE[u(s)]\end{pmatrix},
                                \begin{pmatrix}\dbE[X_0^u(s)]\\\dbE[u(s)]\end{pmatrix}\rran ds\Bigg\}\\
\ns\ds=\blan G\h\cL_tu,\h\cL_tu\bran+\blan \bar G\dbE[\h\cL_tu],\dbE[\h\cL_tu]\bran
+\blan Q\cL_tu,\cL_tu\bran+2\blan S\cL_tu,u\bran+\blan Ru,u\bran\\
\ns\ds\q~+\blan \bar Q\dbE[\cL_tu],\dbE[\cL_tu]\bran+2\blan\bar S\dbE[\cL_tu],\dbE[u]\bran+\blan\bar R\dbE[u],\dbE[u]\bran\\
\ns\ds=\blan\big[\h\cL_t^*(G+\dbE^*\bar G\dbE)\h\cL_t+\cL_t^*(Q+\dbE^*\bar Q\dbE)\cL_t
+(S+\dbE^*\bar S\dbE)\cL_t+\cL_t^*(S^\top+\dbE^*\bar S^\top\dbE)+(R+\dbE^*\bar R\dbE)\big]u,u\bran.
\ea$$
Denote
\bel{cM_t}\cM_t\deq\h\cL_t^*(G+\dbE^*\bar G\dbE)\h\cL_t+\cL_t^*(Q+\dbE^*\bar Q\dbE)\cL_t
+(S+\dbE^*\bar S\dbE)\cL_t+\cL_t^*(S^\top+\dbE^*\bar S^\top\dbE)+(R+\dbE^*\bar R\dbE),\ee
which is a bounded self-adjoint linear operator on $\cU[t,T]$. Then by Proposition \ref{u+v}, the cost
functional $J(t,\xi;u(\cd))$ can be written as
\bel{rep-cost}\ba{ll}
\ds J(t,\xi;u(\cd))=\lan\cM_tu,u\ran+\lan\cD J(t,\xi;0),u\ran+J(t,\xi;0),\\
\ns\ds\qq\qq\qq~\forall (t,\xi)\in[0,T]\times L^2_{\cF_t}(\Om;\dbR^n),\q\forall u(\cd)\in\cU[t,T].
\ea\ee

\ss

Now let us introduce the following conditions.

\ms

{\bf(H3)} The following holds:
\bel{con}J^0(t,0;u(\cd))\ges0,\qq\forall u(\cd)\in\cU[t,T].\ee

{\bf(H4)} There exists a constant $\d>0$ such that
\bel{uni-con}J^0(t,0;u(\cd))\ges\d\,\dbE\int_t^T|u(s)|^2ds,\qq\forall u(\cd)\in\cU[t,T].\ee

\ss

From \rf{rep-cost}, we see that the map $u(\cd)\mapsto J(t,\xi;u(\cd))$ is convex if and only if
\bel{}\cM_t\ges0,\ee
which is also equivalent to (H3), and  $u(\cd)\mapsto J(t,\xi;u(\cd))$ is uniformly convex if and only if
\bel{}\cM_t\ges \d I,\q \hb{for some } \d>0,\ee
which is also equivalent to (H4). The following result tells us that (H3) is necessary for the finiteness
(and open-loop solvability) of Problem (MF-LQ) at $t$, and (H4) is sufficient for the open-loop solvability
of Problem (MF-LQ) at $t$.

\bp{N-S-cond} \sl Let {\rm(H1)--(H2)} hold and $t\in[0,T)$ be given. We have the following:

\ms

{\rm(i)} If Problem {\rm(MF-LQ)} is finite at $t$, then {\rm(H3)} must hold.

\ms

{\rm(ii)} Suppose {\rm(H4)} holds. Then Problem {\rm(MF-LQ)} is uniquely open-loop solvable at $t$,
and the unique optimal control for the initial pair $(t,\xi)$ is given by
\bel{9-13-u*}u^*(\cd)=-{1\over2}\cM_t^{-1}\cD J(t,\xi;0)(\cd).\ee
Moreover,
\bel{}V(t,\xi)= J(t,\xi;0)-{1\over4}\lt|\cM_t^{-{1\over2}}\cD J(t,\xi;0)\rt|^2.\ee
\ep

\it Proof. \rm (i) We prove the result by contradiction. Suppose that $J^0(t,0;u(\cd))<0$
for some $u(\cd)\in\cU[t,T]$. By Proposition \ref{u+v}, we have
$$J(t,\xi;\l u(\cd))=J(t,\xi;0)+\l^2 J^0(t,0;u(\cd))+\l\dbE\int_t^T\lan \cD J(t,\xi;0)(s) ,u(s)\ran ds,\qq\forall \l\in\dbR.$$
Letting $\l\to\i$, we obtain that
$$V(t,\xi)\les\lim_{\l\to\i}J(t,\xi;\l u(\cd))=-\i,$$
which is a contradiction.

\ms

(ii) Suppose {\rm(H4)} holds. Then the operator $\cM_t$ is invertible, and
$$\ba{lll}
\ds J(t,\xi;u(\cd))\4n&=\4n&\ds\lt|\cM_t^{1\over2}u+{1\over2}\cM_t^{-{1\over2}}\cD J(t,\xi;0)\rt|^2
+J(t,\xi;0)-{1\over4}\lt|\cM_t^{-{1\over2}}\cD J(t,\xi;0)\rt|^2,\\
\ns\4n&\ges\4n&\ds J(t,\xi;0)-{1\over4}\lt|\cM_t^{-{1\over2}}\cD J(t,\xi;0)\rt|^2,
\qq~\forall \xi\in L^2_{\cF_t}(\Om;\dbR^n),\q\forall u(\cd)\in\cU[t,T].\ea$$
Note that the equality in the above holds if and only if
$$u=-{1\over2}\cM_t^{-1}\cD J(t,\xi;0).$$
The result therefore follows. \endpf

\ms

Due to the necessity of (H3) for the finiteness of Problem (MF-LQ), we will assume (H3) holds in the rest of this paper.
Now for any $\e>0$, consider state equation \rf{state} and the following cost functional:
\bel{cost-e}\ba{lll}
\ds J_\e(t,\xi;u(\cd))\4n&\deq\4n&\ds J(t,\xi;u(\cd))+\e\dbE\int_t^T|u(s)|^2ds\\
\ns\4n&=\4n&\ds\lan(\cM_t+\e I)u,u\ran+\lan\cD J(t,\xi;0),u\ran+J(t,\xi;0).\ea\ee
Denote the corresponding optimal control problem and value function by Problem (MF-LQ)$_\e$ and $V_\e(\cd\,,\cd)$,
respectively. By Proposition \ref{N-S-cond}, part (ii), for any $\xi\in L^2_{\cF_t}(\Om;\dbR^n)$, Problem (MF-LQ)$_\e$
admits a unique optimal control
\bel{u-e-star}u_\e^*(\cd)=-{1\over2}(\cM_t+\e I)^{-1}\cD J(t,\xi;0)(\cd),\ee
and the value function is given by
\bel{V-e}V_\e(t,\xi)=J(t,\xi;0)-{1\over4}\lt|(\cM_t+\e I)^{-{1\over2}}\cD J(t,\xi;0)\rt|^2.\ee
Now, we are ready to state the main result of this section.

\bt{F-S} \sl Let {\rm(H1)--(H3)} hold and $\xi\in L^2_{\cF_t}(\Om;\dbR^n)$. We have the following:

\ms

{\rm(i)} $\lim_{\e\to 0}V_\e(t,\xi)=V(t,\xi)$. In particular, Problem {\rm(MF-LQ)} is finite at
$(t,\xi)$ if and only if $\{V_\e(t,\xi)\}_{\e>0}$ is bounded from below.

\ms

{\rm(ii)} The sequence $\{u_\e^*(\cd)\}_{\e>0}$ defined by {\rm\rf{u-e-star}} is a minimizing sequence
of $u(\cd)\mapsto J(t,\xi;u(\cd))$:
\bel{Mini-seq}\lim_{\e\to0}J(t,\xi;u^*_\e(\cd))=\inf_{u(\cd)\in\cU[t,T]}J(t,\xi;u(\cd))=V(t,\xi).\ee

\ms

{\rm(iii)} The following statements are equivalent:

\vspace{-3pt}

\begin{enumerate}[\rm(a)]
\setlength{\itemindent}{3.4em}
\setlength{\itemsep}{0pt}
\item Problem {\rm(MF-LQ)} is open-loop solvable at $(t,\xi)$;
\item The sequence $\{u^*_\e(\cd)\}_{\e>0}$ is bounded in $\cU[t,T]$;
\item The sequence $\{u^*_\e(\cd)\}_{\e>0}$ admits a weakly convergent subsequence;
\item The sequence $\{u^*_\e(\cd)\}_{\e>0}$ admits a strongly convergent subsequence.
\end{enumerate}

\vspace{-3pt}

In this case, the weak {\rm(}strong{\rm)} limit of any weakly {\rm(}strongly{\rm)} convergent subsequence of
$\{u^*_\e(\cd)\}_{\e>0}$ is an optimal control of Problem {\rm(MF-LQ)} at $(t,\xi)$.
\et

To prove Theorem \ref{F-S}, we need the following lemma.

\bl{Hilbert} \sl Let $\cH$ be a Hilbert space with norm $|\cd|$ and $\th,\th_n\in\cH$, $n=1,2,\cdots$.

\ms

{\rm(i)} If $\th_n \to \th$ weakly, then $\ds|\th|\les\liminf_{n\to\i}|\th_n|$.

\ms

{\rm(ii)} $\th_n \to \th$ strongly if and only if
$$|\th_n|\to|\th|\qq \hb{and} \qq\th_n\to\th \hb{\q weakly}.$$
\el

\it Proof of Theorem \rm\ref{F-S}. (i) For any $\e_2>\e_1>0$, we have
$$J_{\e_2}(t,\xi;u(\cd))\ges J_{\e_1}(t,\xi;u(\cd))\ges J(t,\xi;u(\cd)),\qq\forall u(\cd)\in\cU[t,T],$$
which implies that
\bel{}V_{\e_2}(t,\xi)\ges V_{\e_1}(t,\xi)\ges V(t,\xi),\qq\forall \e_2>\e_1>0.\ee
Thus, the limit $\lim_{\e\to 0}V_\e(t,\xi)$ exists and
\bel{barV>V}\bar V(t,\xi)\equiv\lim_{\e\to0}V_\e(t,\xi)\ges V(t,\xi).\ee
On the other hand, for any $K,\d>0$, we can find a $u^\d(\cd)\in\cU[t,T]$, such that
$$ V_\e(t,\xi)\les J(t,\xi;u^\d(\cd))+\e\dbE\int_t^T|u^\d(s)|^2ds\les
\max\{V(t,\xi),-K\}+\d+\e\dbE\int_t^T|u^\d(s)|^2ds.$$
Letting $\e\to 0$, we obtain that
$$\bar V(t,\xi)\les \max\{V(t,\xi),-K\}+\d,\qq\forall K,\d>0,$$
from which we see that
\bel{barV<V}\bar V(t,\xi)\les V(t,\xi).\ee
Combining \rf{barV>V}--\rf{barV<V}, we obtain the desired result.

\ms

(ii) If $V(t,\xi)>-\i$, then by (i), we have
$$\ba{ll}
\ds\e\dbE\int_t^T|u^*_\e(s)|^2ds=J_\e(t,\xi;u^*_\e(\cd))-J(t,\xi;u^*_\e(\cd))=V_\e(t,\xi)-J(t,\xi;u^*_\e(\cd))\\
\ns\ds\qq\qq\qq\q~\les V_\e(t,\xi)-V(t,\xi)\to 0\qq \hb{as}\q \e\to 0.\ea$$
Hence,
$$\lim_{\e\to0}J(t,\xi;u^*_\e(\cd))=\lim_{\e\to0}\bigg[V_\e(t,\xi)-\e\dbE\int_t^T|u^*_\e(s)|^2ds\bigg]=V(t,\xi).$$
If $V(t,\xi)=-\i$, then by (i), we have
$$ J(t,\xi;u^*_\e(\cd))\les J_\e(t,\xi;u^*_\e(\cd))=V_\e(t,\xi)\to -\i\qq \hb{as}\q \e\to 0,$$
and \rf{Mini-seq} still holds.

\ms

(iii) (b) $\Ra$ (c) and (d) $\Ra$ (c) are obvious. We next prove (c) $\Ra$ (a).
Let $\{u^*_{\e_k}(\cd)\}_{k\ges1}$ be a weakly convergent subsequence of $\{u^*_\e(\cd)\}_{\e>0}$
with weak limit $u^*(\cd)$. Then $\{u^*_{\e_k}(\cd)\}_{k\ges1}$ is bounded in $\cU[t,T]$.
For any $u(\cd)\in\cU[t,T]$, we have
\bel{c-a}J(t,\xi;u^*_{\e_k}(\cd))+\e_k\dbE\int_t^T|u^*_{\e_k}(s)|^2ds=V_{\e_k}(t,\xi)\les J(t,\xi;u(\cd))+\e_k\dbE\int_t^T|u(s)|^2ds.\ee
Note that $u(\cd)\mapsto J(t,\xi;u(\cd))$ is sequentially weakly lower semi-continuous.
Letting $k\to\i$ in \rf{c-a}, we obtain
$$J(t,\xi;u^*(\cd))\les \liminf_{k\to\i}J(t,\xi;u^*_{\e_k}(\cd))\les J(t,\xi;u(\cd)),\qq\forall u(\cd)\in\cU[t,T].$$
Hence, $u^*(\cd)$ is an optimal control of Problem (MF-LQ) at $(t,\xi)$.
Now it remains to show (a) $\Ra$ (b) and (a) $\Ra$ (d).
Suppose $v^*(\cd)$ is an optimal control of Problem (MF-LQ) at $(t,\xi)$.
Then for any $\e>0$, we have
$$\left\{\2n\ba{ll}
\ds V_\e(t,\xi)=J_\e(t,\xi;u^*_\e(\cd))\ges V(t,\xi)+\e\dbE\int_t^T|u^*_\e(s)|^2ds,\\
\ns\ds V_\e(t,\xi)\les J_\e(t,\xi;v^*(\cd))=V(t,\xi)+\e\dbE\int_t^T|v^*(s)|^2ds,
\ea\right.$$
from which we see that
\bel{a-d-0}\dbE\int_t^T|u^*_\e(s)|^2ds\les{V_\e(t,\xi)-V(t,\xi)\over\e}\les\dbE\int_t^T|v^*(s)|^2ds,\qq\forall \e>0.\ee
Thus, $\{u^*_\e(\cd)\}_{\e>0}$ is bounded in the Hilbert space $\cU[t,T]$
and hence admits a weakly convergent subsequence $\{u^*_{\e_k}(\cd)\}_{k\ges1}$.
Let $u^*(\cd)$ be the weak limit of $\{u^*_{\e_k}(\cd)\}_{k\ges1}$. By the proof of (c) $\Ra$ (a), we see
that $u^*(\cd)$ is also an optimal control of Problem (MF-LQ) at $(t,\xi)$.
Replacing $v^*(\cd)$ with $u^*(\cd)$ in \rf{a-d-0}, we have
\bel{a-d-1}\dbE\int_t^T|u^*_\e(s)|^2ds\les\dbE\int_t^T|u^*(s)|^2ds,\qq\forall \e>0.\ee
Also, by Lemma \ref{Hilbert}, part (i),
\bel{a-d-2}\dbE\int_t^T|u^*(s)|^2ds\les\liminf_{k\to\i}\dbE\int_t^T|u^*_{\e_k}(s)|^2ds.\ee
Combining \rf{a-d-1}--\rf{a-d-2}, we have
$$\dbE\int_t^T|u^*(s)|^2ds=\lim_{k\to\i}\dbE\int_t^T|u^*_{\e_k}(s)|^2ds.$$
Then it follows from Lemma \ref{Hilbert}, part (ii), that $\{u^*_{\e_k}(\cd)\}_{k\ges1}$ converges to $u^*(\cd)$ strongly.
\endpf

\section{Necessary Conditions for the Uniform Convexity and Riccati Equations}

Theorem \ref{F-S} tells us that in order to solve Problem (MF-LQ), we need only solve mean-filed LQ problems with
uniformly convex cost functionals and then pass to the limit. By Proposition \ref{N-S-cond}, under the uniform
convexity condition (H4), the unique optimal control $u^*(\cd)$ for the initial pair $(t,\xi)$ is determined by
\rf{9-13-u*}. However, such a representation is not easy to compute, since $\cM_t^{-1}$ is in an abstract form
and very complicated. Thus, we would like to find some more explicit form of the optimal control. In this section
we shall investigate uniform convexity of the cost functional and show the necessity of solvability of two
Riccati equations for the uniform convexity of the cost functional.

\ms

First, we present the following result concerning the value function of Problem (MF-LQ)$^0$.

\bp{prop-V0} \sl Let {\rm(H1)--(H2)} and {\rm(H4)} hold. Then there exists a constant $\a\in\dbR$ such that
\bel{prop-V0-0}V^0(s,\xi)\ges\a\dbE\big[|\xi|^2\big],
\qq\forall (s,\xi)\in[t,T]\times L^2_{\cF_s}(\Om;\dbR^n)\hb{ with }\dbE[\xi]=0.\ee
\ep

\it Proof. \rm For any $s\in[t,T]$ and any $u(\cd)\in\cU[s,T]$, we define the {\it zero-extension} of $u(\cd)$
as follows:
\bel{ext}[\,0I_{[t,s)}\oplus u(\cd)](r)=\left\{\2n\ba{ll}0,\qq\ r\in[t,s),\\
\ns\ds u(r),\q r\in[s,T].\ea\right.\ee
Then $v(\cd)\equiv0I_{[t,s)}\oplus u(\cd)\in\cU[t,T]$, and due to
the initial state being 0, the solution $X_0^v(\cd)$ of
$$\left\{\2n\ba{ll}
\ds dX_0^v(r)=\big\{A(r)X_0^v(r)+\bar A(r)\dbE[X_0^v(r)]+B(r)v(r)+\bar B(r)\dbE[v(r)]\big\}dr\\
\ns\ds\qq\qq~~+\big\{C(r)X_0^v(r)+\bar C(r)\dbE[X_0^v(r)]+D(r)v(r)+\bar D(r)\dbE[v(r)]\big\}dW(r),\qq r\in[t,T], \\
\ns\ds X_0^v(t)=0,\ea\right.$$
satisfies $X_0^v(r)=0,~r\in[t,s]$. Hence,
\bel{J-s-con}J^0(s,0;u(\cd))=J^0(t,0;0I_{[t,s)}\oplus u(\cd))\ges
\d\,\dbE\int_t^T\big|[0I_{[t,s)}\oplus u(\cd)](r)\big|^2dr=\d\,\dbE\int_s^T|u(r)|^2dr.\ee
Now, let $(X(\cd),Y(\cd),Z(\cd))$ be the solution of the following (decoupled) MF-FBSDE:
\bel{}\left\{\2n\ba{ll}
\ds dX(r)=\big\{AX+\bar A\dbE[X]\big\}dr+\big\{CX+\bar C\dbE[X]\big\}dW(r),\qq r\in[s,T], \\
\ns\ds dY(r)=-\big\{A^\top Y+\bar A^\top\dbE[Y]+C^\top Z+\bar C^\top\dbE[Z]+QX+\bar Q\dbE[X]\big\}dr+ZdW(r), \qq r\in[s,T],\\
\ns\ds X(s)=\xi,\qq Y(T)=GX(T)+\bar G\dbE[X(T)].
\ea\right.\ee
By Proposition \ref{u+v} and \rf{J-s-con}, we have
\bel{4.5}\ba{ll}
\ds J^0(s,\xi;u(\cd))-J^0(s,\xi;0)\\
\ns\ds=J^0(s,0;u(\cd))+2\dbE\int_s^T\blan B^\top Y\1n+\1n\bar B^\top\dbE[Y]\1n+\1nD^\top Z\1n+\1n\bar D^\top\dbE[Z]\1n+\1nSX\1n+\1n\bar S\dbE[X],u\bran\]dr\\
\ns\ds\ges J^0(s,0;u(\cd))-\d\dbE\int_s^T|u(r)|^2dr
-{1\over\d}\dbE\int_s^T\big|B^\top Y\1n+\1n\bar B^\top\dbE[Y]\1n+\1nD^\top Z\1n+\1n\bar D^\top\dbE[Z]\1n+\1nSX\1n+\1n\bar S\dbE[X]\big|^2dr\\
\ns\ds\ges-{1\over\d}\dbE\int_s^T\big|B^\top Y\1n+\1n\bar B^\top\dbE[Y]\1n+\1nD^\top Z\1n+\1n\bar D^\top\dbE[Z]\1n+\1nSX\1n+\1n\bar S\dbE[X]\big|^2dr.
\ea\ee
If $\dbE[\xi]=0$, then $\dbE[X(\cd)]\equiv0$, and one can verify that
\bel{}X(r)=\dbX(r)\dbX(s)^{-1}\xi,\q Y(r)=\dbY(r)\dbX(s)^{-1}\xi,\q Z(r)=\dbZ(r)\dbX(s)^{-1}\xi,\qq r\in[s,T],\ee
where $\dbX(\cd)$ is the solution to the following $\dbR^{n\times n}$-valued SDE:
\bel{dbX}\left\{\2n\ba{ll}
\ds d\dbX(r)=A(r)\dbX(r)dr+C(r)\dbX(r)dW(r),\qq r\in[0,T], \\
\ns\ds \dbX(0)=I,\ea\right.\ee
and $(\dbY(\cd),\dbZ(\cd))$ is the adapted solution to the following $\dbR^{n\times n}$-valued
backward SDE (BSDE, for short):
\bel{dbY-dbZ}\left\{\2n\ba{ll}
\ds d\dbY(r)=-\big[A(r)^\top\dbY(r)+C(r)^\top\dbZ(r)+Q(r)\dbX(r)\big]dr+\dbZ(r)dW(r), \qq r\in[0,T],\\
\ns\ds \dbY(T)=G\dbX(T).\ea\right.\ee
Note that $\dbX(r)\dbX(s)^{-1}$, $\dbY(r)\dbX(s)^{-1}$ and $\dbZ(r)\dbX(s)^{-1}$ are independent of $\cF_s$.
Thus, $\dbE[X(\cd)]=\dbE[Y(\cd)]=\dbE[Z(\cd)]=0$ and (noting \rf{4.5})
$$\ba{ll}
\ds J^0(s,\xi;u(\cd))\ges J^0(s,\xi;0)-{1\over\d}\dbE\int_s^T\big|B(r)^\top Y(r)+D(r)^\top Z(r)+S(r)X(r)\big|^2dr\\
\ns\ds=\dbE\lt\{\lan GX(T),X(T)\ran+\int_s^T\lan Q(r)X(r),X(r)\ran dr\rt\}
-{1\over\d}\dbE\int_s^T\big|B(r)^\top Y(r)+ D(r)^\top Z(r)+S(r)X(r)\big|^2dr\\
\ns\ds=\dbE\lt\{\xi^\top\lt(\big[\dbX(s)^{-1}\big]^\top\dbX(T)^\top G\dbX(T)\dbX(s)^{-1}
+\int_s^T\big[\dbX(s)^{-1}\big]^\top\dbX(r)^\top Q(r)\dbX(r)\dbX(s)^{-1} dr\rt)\xi\rt\}\\
\ns\ds\q\,-{1\over\d}\dbE\int_s^T\xi^\top\big[\dbX(s)^{-1}\big]^\top\big[B(r)^\top\dbY(r)+ D(r)^\top\dbZ(r)+S(r)\dbX(r)\big]^\top\\
\ns\ds\qq\qq\qq\qq\qq\qq\cd\big[B(r)^\top\dbY(r)+ D(r)^\top\dbZ(r)+S(r)\dbX(r)\big]\dbX(s)^{-1}\xi dr\\
\ns\ds=\dbE\Bigg\{\xi^\top\dbE\Bigg(\big[\dbX(s)^{-1}\big]^\top\dbX(T)^\top G\dbX(T)\dbX(s)^{-1}
+\int_s^T\big[\dbX(s)^{-1}\big]^\top\dbX(r)^\top Q(r)\dbX(r)\dbX(s)^{-1}dr\\
\ns\ds\qq\qq\q~-{1\over\d}\int_s^T\big[\dbX(s)^{-1}\big]^\top\big[B(r)^\top\dbY(r)+ D(r)^\top\dbZ(r)+S(r)\dbX(r)\big]^\top\\
\ns\ds\qq\qq\qq\qq\qq\q\cd\big[B(r)^\top\dbY(r)+ D(r)^\top\dbZ(r)+S(r)\dbX(r)\big]\dbX(s)^{-1} dr\Bigg)\xi\Bigg\}\\
\ns\ds\equiv\dbE\big[\xi^\top\dbM(s)\xi\big].\ea$$
Note that $\dbM(\cd):[t,T]\to \dbS^n$ is continuous. The result therefore follows.
\endpf

\ms

Now, let us introduce the following Riccati equation:
\bel{Ric-1}\left\{\2n\ba{ll}
\ds\dot P+PA+A^\top P+C^\top PC+Q\\
\ns\ds\q-\lt(PB+C^\top PD+S^\top\rt)\lt(R+D^\top PD\rt)^{-1}\lt(B^\top P+D^\top PC+S\rt)=0,\qq \ae~s\in[t,T],\\
\ns\ds P(T)=G.\ea\right.\ee
A solution $P(\cd)$ of \rf{Ric-1} is said to be {\it strongly regular} if
\bel{strong-P}R(s)+D(s)^\top P(s)D(s)\ges \d I,\qq\ae~s\in[t,T],\ee
for some $\d>0$. The Riccati equation \rf{Ric-1} is said to be {\it strongly regularly solvable},
if it admits a strongly regular solution. By a standard argument using Gronwall's inequality,
one can show that if the regular solution of \rf{Ric-1} exists, it must be unique.
Compared with the strongly regular solution, the notion of {\it regular solution}, which is closely related to the
closed-loop strategy, was introduced in \cite{Sun-Yong 2014}. The interested reader is referred to \cite{Sun-Li-Yong}
for further information.

\ms

The following result shows that the strongly regular solvability of the Riccati equation \rf{Ric-1}
is necessary for the uniform convexity of the cost functional.

\bt{P-biyao}\sl Let {\rm(H1)--(H2)} and {\rm(H4)} hold.
Then the Riccati equation \rf{Ric-1} is strongly regularly solvable.
\et

To prove the above result, we need the following lemma, whose proof can be found in \cite{Sun-Li-Yong}.

\bl{lmm-Sun}\sl Let {\rm(H1)--(H2)} hold. For any $\Th(\cd)\in L^2(t,T;\dbR^{m\times n})$,
let $P_\Th(\cd)\in C([t,T];\dbS^n)$ be the solution to the following Lyapunov equation:
\bel{Sun-1}\left\{\2n\ba{ll}
\ds{\dot P}_\Th+P_\Th(A+B\Th)+(A+B\Th)^\top P_\Th+(C+D\Th)^\top P_\Th(C+D\Th)\\
\ns\ds\q~+\Th^\top R\Th+S^\top\Th+\Th^\top S+Q=0,\qq\ae~s\in[t,T],\\
\ns\ds P_\Th(T)=G.\ea\right.\ee
If there exists a constant $\b>0$ such that for all $\Th(\cd)\in L^2(t,T;\dbR^{m\times n})$,
\bel{Sun-2}P_\Th(s),~R(s)+D(s)^\top P_\Th(s)D(s)\ges\b I \qq\ae~s\in[t,T],\ee
then the Riccati equation \rf{Ric-1} is strongly regularly solvable.
\el

\it Proof of Theorem {\rm\ref{P-biyao}}. \rm We only need to show that the condition stated in Lemma \ref{lmm-Sun} holds.
To this end, let $\Th(\cd)\in L^2(t,T;\dbR^{m\times n})$ and $P(\cd)\equiv P_\Th(\cd)$ be the corresponding solution of \rf{Sun-1}.
For any deterministic $u(\cd)\in L^2(t,T;\dbR^m)$, let $X^u(\cd)$ be the solution of
\bel{}\left\{\2n\ba{ll}
\ds dX^u(s)=\big[(A+B\Th)X^u+BuW\big]ds+\big[(C+D\Th)X^u+DuW\big]dW,\qq s\in[t,T], \\
\ns\ds X^u(t)=0,\ea\right.\ee
and set
$$v(\cd)\deq \Th(\cd)X^u(\cd)+u(\cd)W(\cd)\in\cU[t,T].$$
Clearly,
\bel{0mean}\dbE[X^u(s)]=0,\q \dbE[v(s)]=0,\qq s\in[t,T].\ee
By the uniqueness of solutions, $X^u(\cd)$ also solves
\bel{}\left\{\2n\ba{ll}
\ds dX^u(s)=\big\{AX^u\1n+\1n\bar A\dbE[X^u]\1n+\1nBv\1n+\1n\bar B\dbE[v]\big\}ds
+\big\{CX^u\1n+\1n\bar C\dbE[X^u]\1n+\1nDv\1n+\1n\bar D\dbE[v]\big\}dW,\q s\in[t,T], \\
\ns\ds X^u(t)=0.\ea\right.\ee
Thus, by applying It\^o's formula to $s\to\lan P(s)X^u(s),X^u(s)\ran$, we have (noting (H4) and \rf{0mean})
$$\ba{ll}
\ds \d\,\dbE\int_t^T|\Th(s)X^u(s)+u(s)W(s)|^2ds=\d\,\dbE\int_t^T|v(s)|^2ds\les J^0(t,0;v(\cd))\\
\ns\ds=\dbE\left\{\lan GX^u(T),X^u(T)\ran
+\int_t^T\[\lan QX^u,X^u\ran+2\lan SX^u,v\ran+\lan Rv,v\ran\] ds\right\}\\
\ns\ds=\dbE\int_t^T\Big\{\blan\dot PX^u,X^u\bran+\blan P\big[(A+B\Th)X^u+BuW\big],X^u\bran+\blan PX^u,(A+B\Th)X^u+BuW\bran\\
\ns\ds\qq\qq~~+\blan P\big[(C+D\Th)X^u+DuW\big],(C+D\Th)X^u+DuW\bran\\
\ns\ds\qq\qq~~+\blan QX^u,X^u\bran+2\blan SX^u,\Th X^u+uW\bran+\blan R(\Th X^u+uW),\Th X^u+uW\bran\Big\}ds\\
\ns\ds=\dbE\int_t^T\Big\{2\,\blan\big[B^\top P+D^\top PC+S+(R+D^\top PD)\Th\big]X^u,uW\bran+\blan (R+D^\top PD)uW,uW\bran\Big\}ds.\ea$$
Hence, for any $u(\cd)\in L^2(t,T;\dbR^m)$, the following holds:
\bel{P>dI-1}\ba{ll}
\ds\dbE\int_t^T\Big\{2\,\blan\big[B^\top P+D^\top PC+S+(R+D^\top PD-\d I)\Th\big]WX^u,u\bran\\
\ns\ds\qq\q~~+W^2\blan (R+D^\top PD-\d I)u,u\bran\Big\}ds=\d\,\dbE\int_t^T|\Th(s)X^u(s)|^2ds\ges0.\ea\ee
Now, applying It\^o's formula again, we have
$$\left\{\2n\ba{ll}
\ds d\,\dbE\big[W(s)X^u(s)\big]=\Big\{\big[A(s)+B(s)\Th(s)\big]\dbE\big[W(s)X^u(s)\big]+sB(s)u(s)\Big\}ds,\qq s\in[t,T], \\
\ns\ds \dbE\big[W(t)X^u(t)\big]=0.\ea\right.$$
Fix any $u_0\in\dbR^m$, take $u(s)=u_0{\bf 1}_{[t^\prime,t^\prime+h]}(s)$, with $t\les t^\prime<t^\prime+h\les T$. Then
$$\dbE\big[W(s)X^u(s)\big]=\left\{\2n\ba{ll}0,\qq\qq\qq\qq\qq\qq\qq\qq s\in[t,t^\prime],\\
\ns\ds\F(s)\int_t^{s\land(t^\prime+h)}\F(r)^{-1}B(r)ru_0dr,\qq s\in[t^\prime,T],\ea\right.$$
where $\F(\cd)$ is the solution of the following $\dbR^{n\times n}$-valued ordinary differential equation (ODE, for short):
$$\left\{\2n\ba{ll}
\ds \dot\F(s)=\big[A(s)+B(s)\Th(s)\big]\F(s),\qq s\in[0,T], \\
\ns\ds \F(0)=I.\ea\right.$$
Consequently, \rf{P>dI-1} becomes
$$\ba{ll}
\ds\int_{t^\prime}^{t^\prime+h}\Big\{2\,\blan\big[B^\top P+D^\top PC+S+(R+D^\top PD-\d I)\Th\big]\F(s)\int_t^s\F(r)^{-1}B(r)ru_0dr,u_0\bran\\
\ns\ds\qq\qq\q~+s\,\blan (R+D^\top PD-\d I)u_0,u_0\bran\Big\}ds\ges0.\ea$$
Dividing both sides by $h$ and letting $h\to 0$, we obtain
$$t^\prime\blan\big[R(t^\prime)+D(t^\prime)^\top P(t^\prime)D(t^\prime)-\d I\big]u_0,u_0\bran\ges 0,
\qq\forall u_0\in\dbR^m,\q\ae~t^\prime\in[t,T],$$
which implies that
\bel{}R(s)+D(s)^\top P(s)D(s)\ges\d I, \qq\ae~s\in[t,T].\ee
Next, for any $(s,x)\in[t,T]\times\dbR^n$, let $X(\cd)$ be the solution of
$$\left\{\2n\ba{ll}
\ds dX(r)=\big[A(r)+B(r)\Th(r)\big]X(r)dr+\big[C(r)+D(r)\Th(r)\big]X(r)dW(r),\qq r\in[s,T], \\
\ns\ds X(s)=W(s)x,\ea\right.$$
and set
$$w(\cd)\deq \Th(\cd)X(\cd)\in\cU[s,T].$$
Similar to the previous argument, by applying It\^o's formula to $r\to\lan P(r)X(r),X(r)\ran$, we can derive that
\bel{}J^0(s,W(s)x;w(\cd))=\dbE\lan P(s)W(s)x,W(s)x\ran=s\,\lan P(s)x,x\ran.\ee
By Proposition \ref{prop-V0}, we have
$$s\,\lan P(s)x,x\ran=J^0(s,W(s)x;w(\cd))\ges\a\dbE\big[|W(s)x|^2\big]=s\a|x|^2,\qq\forall (s,x)\in[t,T]\times\dbR^n,$$
which implies that $ P(s)\ges\a I,\forall s\in[t,T]$. The proof is completed.
\endpf

\ms

From Theorem \ref{P-biyao}, we see that the Riccati equation \rf{Ric-1} is strongly regularly solvable under the
uniform convexity condition (H4). With the strongly regular solution $P(\cd)$ of \rf{Ric-1}, we may further introduce
the following deterministic LQ optimal control problem.

\ms

Consider the state equation
\bel{state-D}\left\{\2n\ba{ll}
\ds \dot y(s)=\big[A(s)+\bar A(s)\big]y(s)+\big[B(s)+\bar B(s)\big]v(s),\qq s\in[t,T], \\
\ns\ds y(t)= x,\ea\right.\ee
and cost functional
\bel{cost-D}\bar J(t,x;v(\cd))\deq \llan(G+\bar G)y(T),y(T)\rran
+\int_t^T\[\llan\Upsilon y,y\rran+2\llan\G y,v\rran+\llan\Si v,v\rran\]ds,\ee
where
\bel{}\left\{\2n\ba{ll}
\ds\Upsilon=Q+\bar Q+(C+\bar C)^\top P(C+\bar C),\\
\ns\ds\G=(D+\bar D)^\top P(C+\bar C)+S+\bar S,\\
\ns\ds\Si=R+\bar R+(D+\bar D)^\top P(D+\bar D).\ea\right.\ee
We pose the following deterministic LQ problem.

\ms

\bf Problem (DLQ). \rm For any given $(t,x)\in[0,T)\times \dbR^n$, find a $v^*(\cd)\in L^2(t,T;\dbR^m)$, such that
\bel{DLQ}\bar J(t,x;v^*(\cd))=\inf_{v(\cd)\in L^2(t,T;\dbR^m)}\bar J(t,x;v(\cd)).\ee

\ms

Note that the Riccati equation associated with Problem (DLQ) is
\bel{Ric-2}\left\{\2n\ba{ll}
\ds\dot\Pi+\Pi(A+\bar A)+(A+\bar A)^\top\Pi+Q+\bar Q+(C+\bar C)^\top P(C+\bar C)\\
\ns\ds\q-\big[\Pi(B+\bar B)+(C+\bar C)^\top P(D+\bar D)+(S+\bar S)^\top\big]\big[R+\bar R+(D+\bar D)^\top P(D+\bar D)\big]^{-1}\\
\ns\ds\qq\,\cd\big[(B+\bar B)^\top\Pi+(D+\bar D)^\top P(C+\bar C)+(S+\bar S)\big]=0,\qq \ae~s\in[t,T],\\
\ns\ds \Pi(T)=G+\bar G.\ea\right.\ee
We have the following result.

\bt{DLQ-tu}\sl Let {\rm(H1)--(H2)} and {\rm(H4)} hold. Then the map $v(\cd)\mapsto \bar J(t,0;v(\cd))$
is uniformly convex, i.e., there exists a $\l>0$ such that
\bel{DLQ-H4}\bar J(t,0;v(\cd))\ges \l\int_t^T|v(s)|^2ds,\qq\forall v(\cd)\in L^2(t,T;\dbR^m).\ee
Consequently, the strongly regular solution $P(\cd)$ of the Riccati equation \rf{Ric-1} satisfies
\bel{}\Si=R+\bar R+(D+\bar D)^\top P(D+\bar D)\gg0,\ee
and the Riccati equation \rf{Ric-2} admits a unique solution $\Pi(\cd)\in C([t,T];\dbS^n)$.
\et

\it Proof. \rm Let $P(\cd)$ be the strongly regular solution of the Riccati equation \rf{Ric-1} and set
$$\Th=-(R+D^\top PD)^{-1}(B^\top P+D^\top PC+S)\in L^2(t,T;\dbR^{m\times n}).$$
We claim that
\bel{J0=barJ}J^0(t,0;\Th(\cd)X(\cd)+v(\cd))=\bar J(t,0;\Th(\cd)y(\cd)+v(\cd)),\qq\forall v(\cd)\in L^2(t,T;\dbR^m).\ee
To prove \rf{J0=barJ}, take any $v(\cd)\in L^2(t,T;\dbR^m)$, let $y(\cd)$ be the solution of
\bel{DLQ-y}\left\{\2n\ba{ll}
\ds \dot y(s)=\big[A(s)+\bar A(s)\big]y(s)+\big[B(s)+\bar B(s)\big]\big[\Th(s)y(s)+v(s)\big],\qq s\in[t,T], \\
\ns\ds y(t)=0,\ea\right.\ee
and $X(\cd)$ be the solution of
\bel{DLQ-X}\left\{\2n\ba{ll}
\ds dX(s)=\Big\{AX+\bar A\dbE[X]+B(\Th X+v)+\bar B\dbE[\Th X+v]\Big\}ds\\
\ns\ds\qq\qq~+\Big\{CX+\bar C\dbE[X]+D(\Th X+v)+\bar D\dbE[\Th X+v]\Big\}dW(s),\qq s\in[t,T], \\
\ns\ds X(t)=0.\ea\right.\ee
Note that $v(\cd)$ is deterministic. Then
$$\left\{\2n\ba{ll}
\ds d\dbE[X(s)]=\Big\{\big(A+\bar A\big)\dbE[X]+\big(B+\bar B\big)\big(\Th\dbE[X]+v\big)\Big\}ds,\qq s\in[t,T], \\
\ns\ds \dbE[X(t)]=0.\ea\right.$$
By the uniqueness of solutions, we see that
\bel{}\dbE[X(s)]=y(s),\qq s\in[t,T].\ee
Now let $z(\cd)=X(\cd)-\dbE[X(\cd)]$. Then
\bel{}\left\{\2n\ba{ll}
\ds dz(s)=(A\1n+\1nB\Th)zds
+\big\{(C\1n+\1nD\Th)z+(C\1n+\1n\bar C)y+(D\1n+\1n\bar D)(\Th y\1n+\1nv)\big\}dW(s),\qq s\in[t,T],\\
\ns\ds z(t)= 0.\ea\right.\ee
Keep in mind that $v(\cd)$ is deterministic and note that
\bel{}\ba{ll}
\ds0=\dot P+P(A+B\Th)+(A+B\Th)^\top P+(C+D\Th)^\top P(C+D\Th)\\
\ns\ds\q~~+\Th^\top R\Th+S^\top\Th+\Th^\top S+Q.\ea\ee
By applying It\^o's formula to $s\mapsto\lan P(s)z(s),z(s)\ran$, we have (also, noting $\dbE[z]=0$)
$$\ba{ll}
\ds J^0(t,0;\Th(\cd)X(\cd)+v(\cd))\\
\ns\ds=\dbE\bigg\{\llan GX(T),X(T)\rran+\llan \bar G\dbE[X(T)],\dbE[X(T)]\rran\\
\ns\ds\qq~+\int_t^T\Big[\llan QX,X\rran+\llan \bar Q\dbE[X],\dbE[X]\rran+2\llan SX,\Th X+v\rran+2\llan \bar S\dbE[X],\dbE[\Th X+v]\rran\\
\ns\ds\qq\qq\qq+\llan R(\Th X+v),(\Th X+v)\rran+\llan \bar R\dbE[\Th X+v],\dbE[\Th X+v]\rran\Big]ds\bigg\}\\
\ns\ds=\dbE\bigg\{\lan Gz(T)z(T)\ran+\int_t^T\[\lan Qz,z\ran+2\lan Sz,\Th z\ran+\lan R\Th z,\Th z\ran\]ds\bigg\}\\
\ns\ds\q+\llan(G+\bar G)y(T),y(T)\rran+\int_t^T\[\llan (Q+\bar Q)y,y\rran+2\llan (S+\bar S)y,\Th y+v\rran\\
\ns\ds\qq\qq\qq\qq\qq\qq\qq\qq+\llan(R+\bar R)(\Th y+v),\Th y+v\rran\]ds\\
\ns\ds=\dbE\int_t^T\Big\{\blan\dot Pz,z\bran+\blan P(A+B\Th)z,z\bran+\blan Pz,(A+B\Th)z\bran\\
\ns\ds\qq\qq~+\blan P\big[(C+D\Th)z+(C+\bar C)y+(D+\bar D)(\Th y+v)\big],\\
\ns\ds\qq\qq\qq\qq~(C+D\Th)z+(C+\bar C)y+(D+\bar D)(\Th y+v)\bran\\
\ns\ds\qq\qq~+\llan\lt(Q+S^\top\Th+\Th^\top S+\Th^\top R\Th\rt)z,z\rran\Big\}ds\\
\ns\ds\q+\llan(G+\bar G)y(T),y(T)\rran+\int_t^T\[\llan (Q+\bar Q)y,y\rran+2\llan (S+\bar S)y,\Th y+v\rran\\
\ns\ds\qq\qq\qq\qq\qq\qq\qq\qq+\llan(R+\bar R)(\Th y+v),\Th y+v\rran\]ds\\
\ns\ds=\int_t^T\llan P\big[(C+\bar C)y+(D+\bar D)(\Th y+v)\big],(C+\bar C)y+(D+\bar D)(\Th y+v)\rran ds\\
\ns\ds\q+\llan(G+\bar G)y(T),y(T)\rran+\int_t^T\[\llan (Q+\bar Q)y,y\rran+2\llan (S+\bar S)y,\Th y+v\rran\\
\ns\ds\qq\qq\qq\qq\qq\qq\qq\qq+\llan(R+\bar R)(\Th y+v),\Th y+v\rran\]ds\\
\ns\ds=\llan(G+\bar G)y(T),y(T)\rran
+\int_t^T\Big\{\llan\big[Q+\bar Q+(C+\bar C)^\top P(C+\bar C)\big]y,y\rran\\
\ns\ds\qq\qq\qq\qq\qq\qq\q\,~+2\llan\big[(D+\bar D)^\top P(C+\bar C)+S+\bar S\big] y,\Th y+v\rran\\
\ns\ds\qq\qq\qq\qq\qq\qq\q\,~+\llan\big[R+\bar R+(D+\bar D)^\top P(D+\bar D)\big](\Th y+v),\Th y+v\rran\Big\}ds\\
\ns\ds=\bar J(t,0;\Th(\cd)y(\cd)+v(\cd)).\ea$$
Thus, \rf{J0=barJ} holds. Consequently, by (H4), we have
$$\ba{ll}
\ds\bar J(t,0;\Th(\cd)y(\cd)+v(\cd))=J^0(t,0;\Th(\cd)X(\cd)+v(\cd))\\
\ns\ds\ges\d\,\dbE\int_t^T|\Th(s)X(s)+v(s)|^2ds\ges\d\int_t^T\big|\dbE[\Th(s)X(s)+v(s)]\big|^2ds\\
\ns\ds=\d\int_t^T|\Th(s)y(s)+v(s)|^2ds,\qq\forall v(\cd)\in L^2(t,T;\dbR^m),\ea$$
which implies the uniform convexity of $v(\cd)\mapsto \bar J(t,0;v(\cd))$.
The rest of the theorem follows now immediately from \cite[Theorem 4.6]{Sun-Li-Yong}.
\endpf

\section{Sufficiency of the Riccati equations}

In the previous section, we proved that the solvability of the Riccati equations \rf{Ric-1} and \rf{Ric-2}
is necessary for the uniform convexity of the cost functional. In this section, we shall show that it is also
sufficient. Moreover, under the uniform convexity condition, the optimal control can be represented explicitly
as a state feedback form via the solutions of the Riccati equations.

\ms

First we need the following lemma.

\bl{lmm-5.1}\sl Let {\rm(H1)--(H2)} hold. For any $u(\cd)\in\cU[t,T]$, let $X_0^u(\cd)$ be the solution of
\bel{X0u-5.1}\left\{\2n\ba{ll}
\ds dX_0^u(s)=\Big\{A(s)X_0^u(s)+\bar A(s)\dbE[X_0^u(s)]+B(s)u(s)+\bar B(s)\dbE[u(s)]\Big\}ds\\
\ns\ds\qq\qq~~+\Big\{C(s)X_0^u(s)+\bar C(s)\dbE[X_0^u(s)]+D(s)u(s)+\bar D(s)\dbE[u(s)]\Big\}dW(s),\qq s\in[t,T], \\
\ns\ds X_0^u(t)=0.\ea\right.\ee
Then for any $\Th(\cd),\bar\Th(\cd)\in L^2(t,T;\dbR^{m\times n})$, there exists a constant $\g>0$ such that
\bel{lmm-5.1-1}\ba{llll}
\ds\dbE\int_t^T\lt|u(s)-\Th(s)\big(X_0^u(s)-\dbE[X_0^u(s)]\big)\rt|^2ds\1n\4n&\ges\4n&\ds\g\,\dbE\int_t^T|u(s)|^2ds,
\q&\forall u(\cd)\in\cU[t,T],\\
\ns\ds\qq\qq~\int_t^T\lt|\dbE[u(s)]-\bar\Th(s)\dbE[X_0^u(s)]\rt|^2ds\1n\4n&\ges\4n&\ds\g\int_t^T\lt|\dbE[u(s)]\rt|^2ds,
\q&\forall u(\cd)\in\cU[t,T].\ea\ee
\el

\it Proof. \rm Let $\Th(\cd)\in L^2(t,T;\dbR^{m\times n})$. Define a bounded linear operator $\cA:\cU[t,T]\to\cU[t,T]$ by
$$\cA u=u-\Th(X_0^u-\dbE[X_0^u]).$$
Then $\cA$ is bijective and its inverse $\cA^{-1}$ is given by
$$\cA^{-1}u=u+\Th\lt(\wt X_0^u-\dbE\big[\wt X_0^u\big]\rt),$$
where $\wt X_0^u(\cd)$ is the solution of
$$\left\{\2n\ba{ll}
\ds d\wt X_0^u(s)=\Big\{(A+B\Th)\wt X_0^u+(\bar A-B\Th)\dbE\big[\wt X_0^u\big]+Bu+\bar B\dbE[u]\Big\}ds\\
\ns\ds\qq\qq\q+\,\Big\{(C+D\Th)\wt X_0^u+(\bar C-D\Th)\dbE\big[\wt X_0^u\big]+Du+D\dbE[u]\Big\}dW(s),\qq s\in[t,T], \\
\ns\ds\wt X_0^u(t)=0.\ea\right.$$
By the bounded inverse theorem, $\cA^{-1}$ is bounded with norm $\|\cA^{-1}\|>0$. Thus,
$$\ba{ll}
\ds\dbE\int_t^T|u(s)|^2ds=\dbE\int_t^T|(\cA^{-1}\cA u)(s)|^2ds
\les\|\cA^{-1}\|\dbE\int_t^T|(\cA u)(s)|^2ds\\
\ns\ds\qq\qq\qq~~\1n=\|\cA^{-1}\|\dbE\int_t^T\lt|u(s)-\Th(s)\big(X_0^u(s)-\dbE[X_0^u(s)]\big)\rt|^2ds,
\qq\forall u(\cd)\in\cU[t,T],\ea$$
which implies the first inequality in \rf{lmm-5.1-1} with $\g=\|\cA^{-1}\|^{-1}$.

\ms

To prove the second, for any $v(\cd)\in L^2(t,T;\dbR^{m\times n})$,
let $y^v(\cd)$ be the solution to the following ODE:
\bel{lmm-5.1-y}\left\{\2n\ba{ll}
\ds \dot y^v(s)=\big[A(s)+\bar A(s)\big]y^v(s)+\big[B(s)+\bar B(s)\big]v(s),\qq s\in[t,T], \\
\ns\ds y^v(t)=0.\ea\right.\ee
For $\bar\Th(\cd)\in L^2(t,T;\dbR^{m\times n})$, we define a bounded linear operator
$\cB:L^2(t,T;\dbR^m)\to L^2(t,T;\dbR^m)$ by
$$\cB v=v-\bar\Th y^v.$$
Similar to the previous argument, one can show that $\cB$ is invertible and
$$\int_t^T\lt|v(s)-\bar\Th(s)y^v(s)\rt|^2ds\ges{1\over\|\cB^{-1}\|}\int_t^T|v(s)|^2ds,
\qq\forall v(\cd)\in L^2(t,T;\dbR^{m\times n}).$$
Observe that $\dbE[X_0^u(\cd)]$ satisfies \rf{lmm-5.1-y} with $v(\cd)=\dbE[u(\cd)]$.
The result therefore follows.
\endpf

\ms

Now we present the main result of this section, which gives a characterization for the uniform convexity
of the cost functional as well as a feedback representation of the optimal control.

\bt{uni-iff} \sl Let {\rm(H1)--(H2)} hold. Then the map $u(\cd)\mapsto J^0(t,0;u(\cd))$ is uniformly convex
if and only if the Riccati equation \rf{Ric-1} admits a strongly regular solution $P(\cd)$ such that
\bel{Si>0}\Si\equiv R+\bar R+(D+\bar D)^\top P(D+\bar D)\gg0,\ee
and the corresponding Riccati equation \rf{Ric-2} admits a solution $\Pi(\cd)$.
In this case, the unique optimal $u^*(\cd)$ of Problem {\rm(MF-LQ)} at $(t,\xi)$ is given by
\bel{u-star}u^*=\Th\big(X^*-\dbE[X^*]\big)+\bar\Th\dbE[X^*]+\f-\dbE[\f]+\bar\f,\ee
where
\bel{Th-f}\left\{\2n\ba{ll}
\ds\Th=-(R+D^\top PD)^{-1}(B^\top P+D^\top PC+S),\\
\ns\ds\bar\Th=-\Si^{-1}\big[(B+\bar B)^\top\Pi+(D+\bar D)^\top P(C+\bar C)+(S+\bar S)\big],\\
\ns\ds\f=-(R+D^\top PD)^{-1}\big[B^\top\eta+D^\top(\z+P\si)+\rho\big],\\
\ns\ds\bar\f=-\Si^{-1}\big\{(B+\bar B)^\top\bar\eta+(D+\bar D)^\top\big(\dbE[\z]+P\dbE[\si]\big)+\dbE[\rho]+\bar\rho\big\},
\ea\right.\ee
with $(\eta(\cd),\z(\cd))$ and $\bar\eta(\cd)$ being the {\rm(}adapted{\rm)} solutions to the following {\rm BSDE}
\bel{eta-zeta}\left\{\2n\ba{ll}
\ds d\eta(s)=-\big[(A+B\Th)^\top\eta+(C+D\Th)^\top\z+(C+D\Th)^\top P\si\\
\ns\ds\qq\qq\q~+\Th^\top\rho+Pb+q\big]ds+\z dW(s),\qq s\in[t,T],\\
\ns\ds\eta(T)=g,\ea\right.\ee
and ordinary differential equation
\bel{bar-eta}\left\{\2n\ba{ll}
\ds\dot{\bar\eta}+\big[(A+\bar A)+(B+\bar B)\bar\Th\big]^\top\bar\eta
+\bar\Th^\top\Big\{(D+\bar D)^\top\big(P\dbE[\si]+\dbE[\z]\big)+\dbE[\rho]+\bar\rho\Big\}\\
\ns\ds~~\1n+(C+\bar C)^\top\big(P\dbE[\si]+\dbE[\z]\big)+\dbE[q]+\bar q+\Pi\dbE[b]=0, \qq\ae~s\in[t,T],\\
\ns\ds\bar\eta(T)=\dbE[g]+\bar g,\ea\right.\ee
respectively, and $X^*(\cd)$ is the solution of the closed-loop system
\bel{closed-sys}\left\{\2n\ba{ll}
\ds dX^*(s)=\Big\{(A+B\Th)\big(X^*\1n-\dbE[X^*]\big)+\big[(A+\bar A)+(B+\bar B)\bar\Th\big]\dbE[X^*]+b\Big\}ds\\
\ns\ds\qq\qq~+\Big\{(C+D\Th)\big(X^*\1n-\dbE[X^*]\big)+\big[(C+\bar C)+(D+\bar D)\bar\Th\big]\dbE[X^*]+\si\Big\}dW(s),\q s\in[t,T], \\
\ns\ds X^*(t)= \xi.\ea\right.\ee
Moreover, the value $V(t,\xi)$ is given by
\bel{V-rep}\ba{ll}
\ds V(t,\xi)
=\dbE\llan P(t)(\xi-\dbE[\xi])+2\eta(t),\xi-\dbE[\xi]\rran+\lan\Pi(t)\dbE[\xi]+2\bar\eta(t),\dbE[\xi]\ran\\
\ns\ds\qq\q\ ~~+\dbE\int_t^T\Big\{\lan P\si,\si\ran+2\lan\eta,b-\dbE[b]\ran+2\lan\z,\si\ran+2\lan\bar\eta,\dbE[b]\ran\\
\ns\ds\qq\qq\qq\qq\q~-\llan\Si_0(\f-\dbE[\f]),\f-\dbE[\f]\rran-\lan\Si\bar\f,\bar\f\ran\Big\}ds,
\ea\ee
where $\Si_0=R+D^\top PD$.
\et

\it Proof. \rm The ``only if\," part has been proved in Section 4. Let us now show the ``if\," part.
For any $\xi\in L^2_{\cF_t}(\Om;\dbR^n)$ and $u(\cd)\in\cU[t,T]$, let $X(\cd)\equiv X(\cd\,;t,\xi,u(\cd))$
be the corresponding solution of \rf{state}. Set
$$z(\cd)=X(\cd)-\dbE[X(\cd)], \q v(\cd)=u(\cd)-\dbE[u(\cd)],\q y(\cd)=\dbE[X(\cd)].$$
Then
\bel{non-z}\left\{\2n\ba{ll}
\ds dz(s)=\Big\{Az+Bv+b-\dbE[b]\Big\}ds\\
\ns\ds\qq\q~~+\Big\{Cz+Dv+\si+(C+\bar C)y+(D+\bar D)\dbE[u]\Big\}dW(s),\qq s\in[t,T], \\
\ns\ds z(t)=\xi-\dbE[\xi],\ea\right.\ee
and
\bel{non-y}\left\{\2n\ba{ll}
\ds \dot y=(A+\bar A)y+(B+\bar B)\dbE[u]+\dbE[b],\qq s\in[t,T],\\
\ns\ds y(t)=\dbE[\xi].\ea\right.\ee
Now we rewrite the cost functional as follows:
%
\bel{J-rewrite}\ba{ll}
\ds J(t,\xi;u(\cd))=\dbE\Bigg\{\lan Gz(T)+2g,z(T)\ran
+\int_t^T\lt[\llan\begin{pmatrix}Q&\1nS^\top\\S&\1nR\end{pmatrix}
                             \begin{pmatrix}z\\v\end{pmatrix},
                             \begin{pmatrix}z\\v\end{pmatrix}\rran
+2\llan\begin{pmatrix}q\\ \rho\end{pmatrix},\begin{pmatrix}z\\v\end{pmatrix}\rran\rt]ds\Bigg\}\\
\ns\ds\qq\qq\qq~\1n+\llan (G+\bar G)y(T)+2\lt(\dbE[g]+\bar g\rt),y(T)\rran\\
\ns\ds\qq\qq\qq~\1n+\int_t^T\lt[\llan\begin{pmatrix}Q\1n+\1n\bar Q&(S\1n+\1n\bar S)^\top\\S\1n+\1n\bar S&R\1n+\1n\bar R\end{pmatrix}
                              \begin{pmatrix}y\\ \dbE[u]\end{pmatrix},
                              \begin{pmatrix}y\\\dbE[u]\end{pmatrix}\rran
+2\llan\begin{pmatrix}\dbE[q]\1n+\1n\bar q\\\dbE[\rho]\1n+\1n\bar\rho\end{pmatrix},
      \begin{pmatrix}y\\\dbE[u]\end{pmatrix}\rran\rt]ds.\ea\ee
Applying It\^o's formula to $s\mapsto\lan P(s)z(s)+2\eta(s),z(s)\ran$, we have (noting $\dbE[z]\equiv0, \dbE[v]\equiv0$)
$$\ba{ll}
\dbE\lan Gz(T)+2g,z(T)\ran-\dbE\lan P(t)(\xi-\dbE[\xi])+2\eta(t),\xi-\dbE[\xi]\ran\\
\ns\ds\q+\,\dbE\int_t^T\[\lan Qz,z\ran+2\lan Sz,v\ran+\lan Rv,v\ran+2\lan q,z\ran+2\lan \rho,v\ran\]ds\\
\ns\ds=\dbE\int_t^T\[\blan\dot Pz,z\bran+\blan P(Az+Bv+b-\dbE[b]),z\bran+\blan Pz,(Az+Bv+b-\dbE[b])\bran\\
\ns\ds\qq\qq~+\blan P\big\{Cz+Dv+\si+(C+\bar C)y+(D+\bar D)\dbE[u]\big\},\\
\ns\ds\qq\qq\qq\qq~Cz+Dv+\si+(C+\bar C)y+(D+\bar D)\dbE[u]\bran\\
\ns\ds\qq\qq~-2\blan(A+B\Th)^\top\eta+(C+D\Th)^\top\z+(C+D\Th)^\top P\si+\Th^\top\rho+Pb+q,z\bran\\
\ns\ds\qq\qq~+2\blan\eta,Az+Bv+b-\dbE[b]\bran\\
\ns\ds\qq\qq~+2\blan\z,Cz+Dv+\si+(C+\bar C)y+(D+\bar D)\dbE[u]\bran\]ds\\
\ns\ds\q+\,\dbE\int_t^T\[\lan Qz,z\ran+2\lan Sz,v\ran+\lan Rv,v\ran+2\lan q,z\ran+2\lan \rho,v\ran\]ds\\
\ns\ds=\dbE\int_t^T\[\blan\big(\dot P+PA+A^\top P+C^\top PC+Q\big)z,z\bran+2\blan\big(PB+C^\top PD+S^\top\big)v,z\bran\\
\ns\ds\qq\qq~+\blan(R+D^\top PD)v,v\bran+2\blan B^\top\eta+D^\top\z+D^\top P\si+\rho,v-\Th z\bran\\
\ns\ds\qq\qq~+2\blan P\dbE[\si]+\dbE[\z],(C+\bar C)y+(D+\bar D)\dbE[u]\bran\\
\ns\ds\qq\qq~+\blan P\big\{(C+\bar C)y+(D+\bar D)\dbE[u]\big\},(C+\bar C)y+(D+\bar D)\dbE[u]\bran\\
\ns\ds\qq\qq~+\lan P\si,\si\ran+2\lan\eta,b-\dbE[b]\ran+2\lan\z,\si\ran\]ds\\
\ns\ds=\dbE\int_t^T\[\blan\Th^\top\Si_0\Th z,z\bran-2\blan\Th^\top\Si_0v,z\bran+\lan\Si_0v,v\ran-2\lan\Si_0\f,v-\Th z\ran\\
\ns\ds\qq\qq~+\blan(C+\bar C)^\top P(C+\bar C)y,y\bran+2\blan(C+\bar C)^\top P(D+\bar D)\dbE[u],y\bran
\ea$$
\bel{Ito-1}\ba{ll}
\ds\qq\qq~+\blan(D+\bar D)^\top P(D+\bar D)\dbE[u],\dbE[u]\bran+2\blan(C+\bar C)^\top\big(P\dbE[\si]+\dbE[\z]\big),y\bran\\
\ns\ds\qq\qq~+2\blan(D+\bar D)^\top\big(P\dbE[\si]+\dbE[\z]\big),\dbE[u]\bran+\lan P\si,\si\ran+2\lan\eta,b-\dbE[b]\ran+2\lan\z,\si\ran\]ds\\
\ns\ds=\dbE\int_t^T\[\lan\Si_0(v-\Th z-\f),v-\Th z-\f\ran-\lan\Si_0\f,\f\ran\\
\ns\ds\qq\qq~+\blan(C+\bar C)^\top P(C+\bar C)y,y\bran+2\blan(C+\bar C)^\top P(D+\bar D)\dbE[u],y\bran\\
\ns\ds\qq\qq~+\blan(D+\bar D)^\top P(D+\bar D)\dbE[u],\dbE[u]\bran+2\blan(C+\bar C)^\top\big(P\dbE[\si]+\dbE[\z]\big),y\bran\\
\ns\ds\qq\qq~+2\blan(D+\bar D)^\top\big(P\dbE[\si]+\dbE[\z]\big),\dbE[u]\bran+\lan P\si,\si\ran+2\lan\eta,b-\dbE[b]\ran+2\lan\z,\si\ran\]ds.
\ea\ee
Applying the integration by parts formula to $s\mapsto\lan\Pi(s)y(s)+2\bar\eta(s),y(s)\ran$, we have
\bel{Ito-2}\ba{ll}
\ds \blan(G+\bar G)y(T)+2(\dbE[g]+\bar g),y(T)\bran-\blan\Pi(t)\dbE[\xi]+2\bar\eta(t),\dbE[\xi]\bran\\
\ns\ds\q\ +\int_t^T\lt[\llan\begin{pmatrix}Q+\bar Q&(S+\bar S)^\top\\S+\bar S&R+\bar R\end{pmatrix}
                              \begin{pmatrix}y\\ \dbE[u]\end{pmatrix},
                              \begin{pmatrix}y\\\dbE[u]\end{pmatrix}\rran
+2\llan\begin{pmatrix}\dbE[q]+\bar q\\\dbE[\rho]+\bar\rho\end{pmatrix},
      \begin{pmatrix}y\\\dbE[u]\end{pmatrix}\rran\rt]ds\\
\ns\ds=\int_t^T\[\blan\dot\Pi y,y\bran+\blan\Pi\big\{(A+\bar A)y+(B+\bar B)\dbE[u]+\dbE[b]\big\},y\bran\\
\ns\ds\qq\q~+\blan\Pi y,(A+\bar A)y+(B+\bar B)\dbE[u]+\dbE[b]\bran\\
\ns\ds\qq\q~+2\blan\dot{\bar\eta},y\bran+2\blan\bar\eta,(A+\bar A)y+(B+\bar B)\dbE[u]+\dbE[b]\bran\]ds\\
\ns\ds\q\ +\int_t^T\[\lan(Q\1n+\1n\bar Q)y,y\ran\1n+\1n2\lan(S\1n+\1n\bar S)y,\dbE[u]\ran\1n+\1n\lan(R\1n+\1n\bar R)\dbE[u],\dbE[u]\ran
\1n+\1n2\lan\dbE[q]\1n+\1n\bar q,y\ran\1n+\1n2\lan\dbE[\rho]\1n+\1n\bar\rho,\dbE[u]\ran\]ds\\
\ns\ds=\int_t^T\Big\{\blan\big[\dot\Pi+\Pi(A+\bar A)+(A+\bar A)^\top\Pi+Q+\bar Q\big]y,y\bran\\
\ns\ds\qq\q~+2\blan\big[\Pi(B+\bar B)+(S+\bar S)^\top\big]\dbE[u],y\bran
+2\blan\dot{\bar\eta}+(A+\bar A)^\top\bar\eta+\dbE[q]+\bar q+\Pi\dbE[b],y\bran\\
\ns\ds\qq\q~+2\blan(B+\bar B)^\top\bar\eta+\dbE[\rho]+\bar\rho,\dbE[u]\bran
+\lan(R+\bar R)\dbE[u],\dbE[u]\ran+2\lan\bar\eta,\dbE[b]\ran\Big\}ds.\ea\ee
Adding \rf{Ito-1} and \rf{Ito-2} together and noting \rf{J-rewrite}, we obtain
$$\ba{ll}
\ds J(t,\xi;u(\cd))-\dbE\lan P(t)(\xi-\dbE[\xi])+2\eta(t),\xi-\dbE[\xi]\ran-\lan\Pi(t)\dbE[\xi]+2\bar\eta(t),\dbE[\xi]\ran\\
\ns\ds=\dbE\int_t^T\Big\{\lan\Sigma_0(v-\Th z-\f),v-\Th z-\f\ran-\lan\Sigma_0\f,\f\ran\\
\ns\ds\qq\qq~+\lan P\si,\si\ran+2\lan\eta,b-\dbE[b]\ran+2\lan\z,\si\ran+2\lan\bar\eta,\dbE[b]\ran\Big\}ds\\
\ns\ds\q+\int_t^T\Big\{\blan\big[\dot\Pi+\Pi(A+\bar A)+(A+\bar A)^\top\Pi+Q+\bar Q+(C+\bar C)^\top P(C+\bar C)\big]y,y\bran\\
\ns\ds\qq\qq~+2\blan\big[\Pi(B+\bar B)+(C+\bar C)^\top P(D+\bar D)+(S+\bar S)^\top\big]\dbE[u],y\bran\\
\ns\ds\qq\qq~+\blan\big[R+\bar R+(D+\bar D)^\top P(D+\bar D)\big]\dbE[u],\dbE[u]\bran\\
\ns\ds\qq\qq~+2\blan\dot{\bar\eta}+(A+\bar A)^\top\bar\eta+(C+\bar C)^\top\big(P\dbE[\si]+\dbE[\z]\big)+\dbE[q]+\bar q+\Pi\dbE[b],y\bran\\
\ns\ds\qq\qq~+2\blan(B+\bar B)^\top\bar\eta+(D+\bar D)^\top\big(P\dbE[\si]+\dbE[\z]\big)+\dbE[\rho]+\bar\rho,\dbE[u]\bran\Big\}ds\\
\ns\ds=\dbE\int_t^T\Big\{\lan\Sigma_0(v-\Th z-\f),v-\Th z-\f\ran-\lan\Sigma_0\f,\f\ran\\
\ns\ds\qq\qq~+\lan P\si,\si\ran+2\lan\eta,b-\dbE[b]\ran+2\lan\z,\si\ran+2\lan\bar\eta,\dbE[b]\ran\Big\}ds\\
\ns\ds\q+\int_t^T\Big\{\blan\bar\Th^\top\Si\bar\Th y,y\bran-2\blan\bar\Th^\top\Si\dbE[u],y\bran
+\lan\Si\dbE[u],\dbE[u]\ran-2\lan\Si\bar\f,\dbE[u]-\bar\Th y\ran\Big\}ds\\
\ns\ds=\dbE\int_t^T\Big\{\lan\Sigma_0(v-\Th z-\f),v-\Th z-\f\ran-\lan\Sigma_0\f,\f\ran
\ea$$
\bel{Ito-1+2}\ba{ll}
\ds\qq\qq~+\blan\Si\big(\dbE[u]-\bar\Th y-\bar\f\big),\dbE[u]-\bar\Th y-\bar\f\bran-\lan\Si\bar\f,\bar\f\ran\\
\ns\ds\qq\qq~+\lan P\si,\si\ran+2\lan\eta,b-\dbE[b]\ran+2\lan\z,\si\ran+2\lan\bar\eta,\dbE[b]\ran\Big\}ds\\
\ns\ds=\dbE\int_t^T\Big\{\lan P\si,\si\ran+2\lan\eta,b-\dbE[b]\ran
+2\lan\z,\si\ran+2\lan\bar\eta,\dbE[b]\ran-\lan\Si\bar\f,\bar\f\ran\\
\ns\ds\qq\qq~-\lan\Si_0(\f-\dbE[\f]),\f-\dbE[\f]\ran+\lan\Si_0(v-\Th z-\f+\dbE[\f]),v-\Th z-\f+\dbE[\f]\ran\\
\ns\ds\qq\qq~+\blan\Si\big(\dbE[u]-\bar\Th y-\bar\f\big),\dbE[u]-\bar\Th y-\bar\f\bran\Big\}ds.
\ea\ee
Since $\Si_0,\Si\gg0$, \rf{Ito-1+2} implies that
\bel{Ito-1+2=3}\ba{ll}
\ds J(t,\xi;u(\cd))\ges\dbE\lan P(t)(\xi-\dbE[\xi])+2\eta(t),\xi-\dbE[\xi]\ran+\lan\Pi(t)\dbE[\xi]+2\bar\eta(t),\dbE[\xi]\ran\\
\ns\ds\qq\qq\qq+\,\dbE\int_t^T\Big\{\lan P\si,\si\ran+2\lan\eta,b-\dbE[b]\ran+2\lan\z,\si\ran+2\lan\bar\eta,\dbE[b]\ran\\
\ns\ds\qq\qq\qq\qq\qq~-\lan\Si_0(\f-\dbE[\f]),\f-\dbE[\f]\ran-\lan\Si\bar\f,\bar\f\ran\Big\}ds,
\ea\ee
with the equality holding if and only if
$$\left\{\2n\ba{ll}
\ds u-\dbE[u]=v=\Th z+\f-\dbE[\f]=\Th\big(X-\dbE[X]\big)+\f-\dbE[\f],\\
\ns\ds \dbE[u]=\bar\Th y+\bar\f=\bar\Th\dbE[X]+\bar\f,\ea\right.$$
which is also equivalent to
\bel{}u=\Th\big(X-\dbE[X]\big)+\bar\Th\dbE[X]+\f-\dbE[\f]+\bar\f.\ee
In particular, when $b(\cd),\si(\cd),g(\cd),\bar g(\cd),q(\cd), \bar q(\cd),\rho(\cd),\bar \rho(\cd)=0$, we have
$$(\eta(\cd),\z(\cd))=(0,0),\q\bar\eta(\cd)=0,\q \f(\cd)=\bar\f(\cd)=0.$$
Take $\xi=0$. Then $X(\cd)$ satisfies
\bel{}\left\{\2n\ba{ll}
\ds dX(s)=\Big\{A(s)X(s)+\bar A(s)\dbE[X(s)]+B(s)u(s)+\bar B(s)\dbE[u(s)]\Big\}ds\\
\ns\ds\qq\qq~~+\Big\{C(s)X(s)+\bar C(s)\dbE[X(s)]+D(s)u(s)+\bar D(s)\dbE[u(s)]\Big\}dW(s),\qq s\in[t,T], \\
\ns\ds X(t)=0,\ea\right.\ee
and \rf{Ito-1+2} becomes
\bel{Ito-J0}\ba{ll}
\ds J^0(t,0;u(\cd))=\dbE\int_t^T\Big\{\blan\Si_0\big[u-\dbE[u]-\Th(X-\dbE[X])\big],u-\dbE[u]-\Th(X-\dbE[X])\bran\\
\ns\ds\qq\qq\qq\qq\qq~+\blan\Si\big(\dbE[u]-\bar\Th\dbE[X]\big),\dbE[u]-\bar\Th\dbE[X]\bran\Big\}ds.
\ea\ee
Noting that $\Si_0,\Si\ges\d I$ for some $\d>0$ and making use of Lemma \ref{lmm-5.1}, we have
\bel{}\ba{ll}
\ds J^0(t,0;u(\cd))\ges\d\,\dbE\int_t^T\Big\{\lt|u-\dbE[u]-\Th(X-\dbE[X])\rt|^2+\lt|\dbE[u]-\bar\Th\dbE[X]\rt|^2\Big\}ds\\
\ns\ds\qq\qq\q~\ges\d\,\dbE\int_t^T|u-\Th(X-\dbE[X])|^2-2\lan u-\Th(X-\dbE[X]),\dbE[u]\ran+(1+\g)|\dbE[u]|^2 ds\\
\ns\ds\qq\qq\q~\ges{\d\g\over1+\g}\,\dbE\int_t^T|u-\Th(X-\dbE[X])|^2 ds
\ges{\d\g^2\over1+\g}\,\dbE\int_t^T|u(s)|^2 ds,\qq\forall u(\cd)\in\cU[t,T],\ea\ee
for some $\g>0$. The uniform convexity of $u(\cd)\mapsto J^0(t,0;u(\cd))$ follows immediately.
\endpf

\ms

Note that for Problem (MF-LQ)$^0$ (where $b(\cd),\si(\cd),g(\cd),\bar g(\cd),q(\cd), \bar q(\cd),\rho(\cd),\bar \rho(\cd)=0$),
under the uniform convexity condition (H4), the value at $(t,\xi)$ is given by
\bel{} V^0(t,\xi)=\dbE\llan P(t)(\xi-\dbE[\xi]),\xi-\dbE[\xi]\rran+\lan\Pi(t)\dbE[\xi],\dbE[\xi]\ran,\ee
where $P(\cd)$ and $\Pi(\cd)$ are the solutions to the Riccati equations \rf{Ric-1} and \rf{Ric-2}, respectively.
The unique optimal $u^*(\cd)$ is given by
\bel{}u^*=\Th\big(X^*-\dbE[X^*]\big)+\bar\Th\dbE[X^*],\ee
where $\Th(\cd), \bar\Th(\cd)$ are defined by \rf{Th-f} and $X^*(\cd)$ is the solution of
\bel{}\left\{\2n\ba{ll}
\ds dX^*(s)=\Big\{(A+B\Th)\big(X^*\1n-\dbE[X^*]\big)+\big[(A+\bar A)+(B+\bar B)\bar\Th\big]\dbE[X^*]\Big\}ds\\
\ns\ds\qq\qq~+\Big\{(C+D\Th)\big(X^*\1n-\dbE[X^*]\big)+\big[(C+\bar C)+(D+\bar D)\bar\Th\big]\dbE[X^*]\Big\}dW(s),\q s\in[t,T], \\
\ns\ds X^*(t)= \xi.\ea\right.\ee

\ms

To conclude this section, we present a sufficient condition for the uniform convexity of the cost functional.
From the following result, we will see that \rf{Classic} implies the uniform convexity condition (H4). However,
the converse fails. A counterexample will be present in the next section (see Example \ref{exm1}).

\bp{suf-con} \sl Let {\rm(H1)--(H2)} hold and $t\in[0,T)$ be given. If there exists a constant $\d>0$ such that
\bel{standard}\left\{\2n\ba{ll}
\ds G,\,G+\bar G\ges0,\q R(s),\,R(s)+\bar R(s)\ges\d I,\q Q(s)-S(s)^\top R(s)^{-1}S(s)\ges0,\\
\ns\ds Q(s)+\bar Q(s)-\big[S(s)+\bar S(s)\big]^\top\big[R(s)+\bar R(s)\big]^{-1}\big[S(s)+\bar S(s)\big]\ges0,
\ea\right.\q\ae~s\in[t,T],\ee
then the map $u(\cd)\mapsto J^0(t,0;u(\cd))$ is uniformly convex.
\ep

\it Proof. \rm For any $u(\cd)\in\cU[t,T]$, let $X_0^u(\cd)$ be the solution of \rf{X0u-5.1}. Then
$$\ba{ll}
\ds J^0(t,0;u(\cd))=\dbE\Bigg\{\llan G\big(X_0^u(T)-\dbE[X_0^u(T)]\big),X_0^u(T)-\dbE[X_0^u(T)]\rran\\
\ns\ds\qq\qq\qq\qq~+\int_t^T\llan\begin{pmatrix}Q&S^\top\\S&R\end{pmatrix}
                             \begin{pmatrix}X_0^u-\dbE[X_0^u]\\u-\dbE[u]\end{pmatrix},
                             \begin{pmatrix}X_0^u-\dbE[X_0^u]\\u-\dbE[u]\end{pmatrix}\rran ds\Bigg\}\\
\ns\ds\qq\qq\qq~+\llan (G+\bar G)\dbE[X_0^u(T)],\dbE[X_0^u(T)]\rran\\
\ns\ds\qq\qq\qq\qq~+\int_t^T\llan\begin{pmatrix}Q+\bar Q&(S+\bar S)^\top\\S+\bar S&R+\bar R\end{pmatrix}
                              \begin{pmatrix}\dbE[X_0^u]\\ \dbE[u]\end{pmatrix},
                              \begin{pmatrix}\dbE[X_0^u]\\\dbE[u]\end{pmatrix}\rran ds\\
\ns\ds\ges\dbE\int_t^T\Big\{\llan Q\big(X_0^u-\dbE[X_0^u]\big),X_0^u-\dbE[X_0^u]\rran+2\llan S\big(X_0^u-\dbE[X_0^u]\big),u-\dbE[u]\rran\\
\ns\ds\qq\qq~~+\llan R\big(u-\dbE[u]\big),u-\dbE[u]\rran\Big\} ds\\
\ns\ds~~~+\int_t^T\Big\{\llan \big(Q+\bar Q\big)\dbE[X_0^u],\dbE[X_0^u]\rran
+2\llan\big(S+\bar S\big)\dbE[X_0^u],\dbE[u]\rran+\llan\big(R+\bar R\big)\dbE[u],\dbE[u]\rran\Big\} ds\\
\ns\ds=\dbE\int_t^T\Big\{\llan\big(Q-S^\top R^{-1}S\big)\big(X_0^u-\dbE[X_0^u]\big),X_0^u-\dbE[X_0^u]\rran\\
\ns\ds\qq\qq~~+\llan R\big[u-\dbE[u]+R^{-1}S\big(X_0^u-\dbE[X_0^u]\big)\big],
u-\dbE[u]+R^{-1}S\big(X_0^u-\dbE[X_0^u]\big)\rran\Big\}ds\\
\ns\ds~~~+\int_t^T\Big\{\llan\big[Q+\bar Q-(S+\bar S)^\top(R+\bar R)^{-1}(S+\bar S)\big]\dbE[X_0^u],\dbE[X_0^u]\rran\\
\ns\ds\qq\qq~~+\llan\big(R+\bar R\big)\(\dbE[u]+(R+\bar R)^{-1}(S+\bar S)\dbE[X_0^u]\),\dbE[u]+(R+\bar R)^{-1}(S+\bar S)\dbE[X_0^u]\rran\Big\}ds\\
\ns\ds\ges\d\,\dbE\int_t^T\lt|u-\dbE[u]+R^{-1}S\big(X_0^u-\dbE[X_0^u]\big)\rt|^2ds
+\d\int_t^T\lt|\dbE[u]+(R+\bar R)^{-1}(S+\bar S)\dbE[X_0^u]\rt|^2ds.\ea$$
Consequently, by Lemma \ref{lmm-5.1} (taking $\Th=-R^{-1}S$ and $\bar\Th=-(R+\bar R)^{-1}(S+\bar S)$), we have
\bel{}\ba{ll}
\ds J^0(t,0;u(\cd))\ges\d\,\dbE\int_t^T\Big\{\lt|u-\dbE[u]+R^{-1}S\big(X_0^u-\dbE[X_0^u]\big)\rt|^2+\g|\dbE[u]|^2\Big\}ds\\
\ns\ds\qq\qq\q~\ges{\d\g\over1+\g}\dbE\int_t^T\lt|u+R^{-1}S\big(X_0^u-\dbE[X_0^u]\big)\rt|^2ds\\
\ns\ds\qq\qq\q~\ges{\d\g^2\over1+\g}\dbE\int_t^T|u(s)|^2ds,\qq\forall u(\cd)\in\cU[t,T],\ea\ee
for some $\g>0$. This completes the proof.
\endpf

\section{Examples}

In this section we present two illustrative examples. In the first example, the condition \rf{Classic} does not hold,
but the corresponding Riccati equations are still solvable. Thus, by Theorem \ref{uni-iff}, the cost functional is uniformly
convex. This example shows that the uniform convexity condition (H4) is indeed weaker than \rf{Classic}.

\bex{exm1}\rm  Consider the following Problem (MF-LQ)$^0$ with one-dimensional state equation
\bel{}\left\{\2n\ba{ll}
\ds dX(s)=\big\{\dbE[X(s)]+u(s)+\dbE[u(s)]\big\}ds+\sqrt{2}u(s)dW(s),\qq s\in[t,1], \\
\ns\ds X(t)=\xi,\ea\right.\ee
and cost functional
\bel{}J(t,\xi;u(\cd))=\dbE\lt\{G|X(1)|^2+\bar G\lt|\dbE[X(1)]\rt|^2+\int_t^1\(R(s)|u(s)|^2+\bar R(s)|\dbE[u(s)]|^2\)ds\rt\},\ee
where
$$\left\{\2n\ba{ll}
G=8,\qq \bar G=-\a-8 \q \hb{with}\q 0<\a<{1\over 2(e^2-1)},\\
\ns R(s)=(s+1)^3-4(s+1)^2,\qq \bar R(s)=1-(s+1)^3,\qq s\in[0,1].
\ea\right.$$
The Riccati equations for the above problem are
$$\left\{\2n\ba{ll}
\ds \dot P(s)-{P(s)^2\over R(s)+2P(s)}=0, \\
\ns\ds P(1)=8,\ea\right.$$
and
$$\left\{\2n\ba{ll}
\ds \dot\Pi(s)+2\Pi(s)-{4\Pi(s)^2\over R(s)+\bar R(s)+2P(s)}=0, \\
\ns\ds \Pi(1)=-\a.\ea\right.$$
Clearly,
$$G+\bar G=-\a<0,\q R(s)=(s+1)^2(s-3)\les-2,\q R(s)+\bar R(s)=1-4(s+1)^2\les-3,\qq s\in[0,1].$$
Hence, the condition \rf{Classic} does not hold. However, one can verify that the above Riccati equations
are solvable on the whole interval $[0,1]$ with solutions given by
$$P(s)=2(s+1)^2,\qq \Pi(s)={\a e^{2(1-s)}\over 2\a[e^{2(1-s)}-1]-1}<0,\qq s\in[0,1].$$
Note that
$$\left\{\2n\ba{ll}
\ds R(s)+D(s)^\top P(s) D(s)=(s+1)^3\ges 1, \\
\ns\ds R(s)+\bar R(s)+\big[D(s)+\bar D(s)\big]^\top P(s)\big[D(s)+\bar D(s)\big]=1,
\ea\right. \qq s\in[0,1].$$
By Theorem \ref{uni-iff}, the cost functional $J(t,\xi;u(\cd))$ is uniformly convex in $u(\cd)$,
and for any initial pair $(t,\xi)\in[0,1)\times L^2_{\cF_t}(\Om;\dbR)$, the problem admits a
unique optimal control $u^*(\cd)$ given by
$$u^*(s)=-{2\over s+1}X^*(s)+\lt[{2\over s+1}-2\Pi(s)\rt]\dbE[X^*(s)],\qq s\in[t,1]$$
with $X^*(\cd)$ being the solution to the following closed-loop system:
$$\left\{\2n\ba{ll}
\ds dX^*(s)=\lt\{-{2\over s+1}X^*(s)+\lt[{s+3\over s+2}-4\Pi(s)\rt]\dbE[X^*(s)]\rt\}ds\\
\ns\ds\qq\qq\q+\lt\{-{2\sqrt{2}\over s+2}X^*(s)+2\sqrt{2}\lt[{1\over s+1}-\Pi(s)\rt]\dbE[X^*(s)]\rt\}dW(s),\qq s\in[t,1], \\
\ns\ds X^*(t)=\xi.\ea\right.$$
\ex

Now we present an example in which the mean-field LQ problem is open-loop solvable,
but the cost functional is not uniformly convex. Hence, the optimal control cannot be constructed directly
in terms of the Riccati equations. However, an optimal could still be found by making use of Theorem
\ref{F-S} and \ref{uni-iff}.

\bex{ex2}\rm Consider the following Problem (MF-LQ)$^0$ with one-dimensional state equation
\bel{}\left\{\2n\ba{ll}
\ds dX(s)=\big\{X(s)-\dbE[X(s)]+\dbE[u(s)]\big\}ds+u(s)dW(s),\qq s\in[t,T], \\
\ns\ds X(t)=\xi,\ea\right.\ee
and cost functional
\bel{}J(t,\xi;u(\cd))=\dbE\lt\{2|X(T)|^2\1n+\1n\lt|\dbE[X(T)]\rt|^2
\1n+\1n\int_t^T\1n\(\1n-4|X(s)|^2\1n-\1n|u(s)|^2\1n+\1n4\lt|\dbE[X(s)]\rt|^2\1n-\1n|\dbE[u(s)]|^2\)ds\rt\}.\ee
In this example,
$$\left\{\2n\ba{ll}
\ds A=1,\q \bar A=-1,\q B=0,\q \bar B=1,\q C=\bar C=0,\q D=1,\q \bar D=0, \\
\ns\ds G=2,\q \bar G=1,\q Q=-4,\q \bar Q=4,\q S=\bar S=0,\q R=\bar R=-1.
\ea\right.$$
Clearly, the condition \rf{Classic} does not hold. The Riccati equations for the problem are
\bel{9-14-P}\left\{\2n\ba{ll}
\ds \dot P(s)+2P(s)-4=0,\qq s\in[t,T], \\
\ns\ds P(T)=2,\ea\right.\ee
and
\bel{9-14-Pi}\left\{\2n\ba{ll}
\ds \dot\Pi(s)-{\Pi(s)^2\over P(s)-2}=0,\qq s\in[t,T], \\
\ns\ds \Pi(T)=3.\ea\right.\ee
It is easy to see that $P(\cd)\equiv2$ is the unique solution of \rf{9-14-P}. However, since
$$R(s)+\bar R(s)+[D(s)+\bar D(s)]^\top P(s)[D(s)+\bar D(s)]=P(s)-2=0,\qq s\in[t,T],$$
we cannot use \rf{9-14-Pi} to solve the problem directly. To investigate the open-loop solvability of the above problem,
let us now consider the following cost functionals for $\e>0$:
\bel{}\ba{ll}
\ds J_\e(t,\xi;u(\cd))=J(t,\xi;u(\cd))+\e\dbE\int_t^T|u(s)|^2ds\\
\ns\ds=\dbE\lt\{2|X(T)|^2+\lt|\dbE[X(T)]\rt|^2
+\int_t^T\(\1n-4|X(s)|^2+(\e-1)|u(s)|^2+4\lt|\dbE[X(s)]\rt|^2-|\dbE[u(s)]|^2\)ds\rt\}.
\ea\ee
We denote the corresponding mean-field LQ problem and value function by Problem (MF-LQ)$_\e^0$
and $V^0_\e(\cd\,,\cd)$, respectively. The Riccati equations for Problem (MF-LQ)$_\e^0$ are
$$\left\{\2n\ba{ll}
\ds \dot P_\e(s)+2P_\e(s)-4=0,\qq s\in[t,T], \\
\ns\ds P_\e(T)=2,\ea\right.$$
and
$$\left\{\2n\ba{ll}
\ds \dot\Pi_\e(s)-{\Pi_\e(s)^2\over\e-2+P_\e(s)}=0,\qq s\in[t,T], \\
\ns\ds \Pi_\e(T)=3.\ea\right.$$
A straightforward calculation leads to
$$P_\e(s)=2,\qq \Pi_\e(s)={3\e\over\e+3(T-s)};\qq s\in[t,T].$$
Since
$$R+\e+D^\top P_\e D=1+\e, \qq R+\e+\bar R+(D+\bar D)^\top P_\e(D+\bar D)=\e,$$
by Theorem \ref{uni-iff}, the map $u(\cd)\mapsto J_\e(t,0;u(\cd))$ is uniformly convex for all $\e>0$
and hence $u(\cd)\mapsto J(t,0;u(\cd))$ is convex. Moreover,
\bel{ex2-Ve0}V_\e^0(t,\xi)=\dbE\llan P_\e(t)(\xi-\dbE[\xi]),\xi-\dbE[\xi]\rran+\lan\Pi_\e(t)\dbE[\xi],\dbE[\xi]\ran,\ee
and the unique optimal control of Problem (MF-LQ)$_\e^0$ at $(t,\xi)$ is given by
\bel{}u_\e^*(s)=-{\Pi_\e(s)\over\e}\dbE[X_\e^*(s)],\qq s\in[t,T],\ee
with $X_\e^*(\cd)$ being the solution to the following closed-loop system:
$$\left\{\2n\ba{ll}
\ds dX_\e^*(s)=\lt\{X_\e^*(s)-\lt(1+{\Pi_\e(s)\over\e}\rt)\dbE[X_\e^*(s)]\rt\}ds
-{\Pi_\e(s)\over\e}\dbE[X_\e^*(s)]dW(s),\qq s\in[t,T], \\
\ns\ds X_\e^*(t)=\xi.\ea\right.$$
Letting $\e\to 0$ in \rf{ex2-Ve0}, we have from Theorem \ref{F-S} that
\bel{}V^0(t,\xi)=\lim_{\e\to0}V_\e^0(t,\xi)=\left\{\2n\ba{ll}
\ds 2\var[\xi],\q &0\les t<T,\\
\ns\ds 2\dbE[\xi^2]+(\dbE[\xi])^2,\q &t=T.\ea\right.\ee
Note that
$$\left\{\2n\ba{ll}
\ds d\dbE[X_\e^*(s)]=-{\Pi_\e(s)\over\e}\dbE[X_\e^*(s)]ds,\qq s\in[t,T], \\
\ns\ds \dbE[X_\e^*(t)]=\dbE[\xi].\ea\right.$$
Hence,
$$\dbE[X_\e^*(s)]=\dbE[\xi]\exp\lt\{-\int_t^s{\Pi_\e(r)\over\e}dr\rt\}
={\e+3(T-s)\over\e+3(T-t)}\dbE[\xi],\qq s\in[t,T],$$
and
$$ u_\e^*(s)=-{\Pi_\e(s)\over\e}\dbE[X_\e^*(s)]=-{3\dbE[\xi]\over\e+3(T-t)}\qq s\in[t,T].$$
It is clear that for $t\in[0,T)$, $u_\e^*(s)$ converges uniformly to
\bel{}u^*(s)\equiv-{\dbE[\xi]\over T-t}\qq s\in[t,T],\ee
which, by Theorem \ref{F-S}, is an optimal control of the original problem at $(t,\xi)$.
\ex

\ms

\bf Acknowledgements. \rm The author wishes to thank Prof. Jiongmin Yong for his valuable comments, 
which have helped to improve the quality of the manuscript. The author also would like to thank 
Dr. Xun Li for his useful suggestions and  financial support.


\end{document}